\documentclass[11pt,a4paper]{article}
\pdfoutput=1
\usepackage[utf8]{inputenc}
\usepackage{amsmath}
\usepackage{amsfonts}
\usepackage{amssymb}
\usepackage{graphicx}
\usepackage{xcolor}
\usepackage[normalem]{ulem}
\usepackage[hyphens]{url}
\usepackage{natbib}
\usepackage{booktabs}
\usepackage{tabularx}
\usepackage{eurosym}
\usepackage{lscape}
\usepackage{soul}

\usepackage{tikz,pgfplots}

\usepackage[font=small]{caption} 
\newcommand{\shortname}{BLTP-RN}
\newcommand{\name}{Bi-objective Long-haul Transportation Problem on a Road Network}
\newcommand{\hos}{HoS}
\newcommand{\gis}{GIS}

\usepackage{algorithmic}
\usepackage{algorithm}
\usepackage[symbol*]{footmisc}

\title{The \name}
\author{Claudia Archetti$^1$, Andrea Mor$^2$, Ola Jabali$^3$\\ Alberto Simonetto$^4$, M. Grazia
Speranza$^2$\\
\vspace{2mm} \\
{\small $^1$ESSEC Business School in Paris,} \\
{\small Department of Information Systems, Decision Sciences and Statistics} \\
{\small Paris, France} \\
{\small archetti@essec.edu}
\vspace{2mm} \\
{\small $^2$University of Brescia,} \\
{\small Department of Economics and Management,} \\
{\small Brescia, Italy} \\
{\small \{andrea.mor,grazia.speranza\}@unibs.it}
\vspace{2mm} \\
{\small $^3$ Politecnico di Milano,}\\
{\small Department of Electronics, Information and Bioengineering}\\
{\small Milan, Italy} \\
{\small ola.jabali@polimi.it}
\vspace{2mm} \\
{\small $^4$ Multiprotexion srl,}\\
{\small Gropello Cairoli, Italy} \\
{\small simonetto@multiprotexion.eu}
}

\begin{document}

\maketitle

\begin{abstract}

In this paper we study a long-haul truck scheduling problem where a path has to be determined for a vehicle traveling from a specified origin to a specified destination. We consider refueling decisions along the path while accounting for heterogeneous fuel prices in a road network. Furthermore, the path has to comply with Hours of Service (\hos{}) regulations. Therefore, a path is defined by the actual road trajectory traveled by the vehicle, as well as the locations where the vehicle stops due to refueling, compliance with \hos{} regulations, or a combination of the two. This setting is cast in a bi-objective optimization problem, considering the minimization of fuel cost and the minimization of path duration. An algorithm is proposed to solve the problem on a road network. The algorithm builds a set of non-dominated paths with respect to the two objectives. Given the enormous theoretical size of the road network, the algorithm follows an interactive path construction mechanism. Specifically, the algorithm dynamically interacts with a geographic information system to identify the relevant potential paths and stop locations.
Computational tests are made on real-sized instances where the distance covered ranges from 500 to 1500 km. The algorithm is compared with solutions obtained from a policy mimicking the current practice of a logistics company. The results show that the non-dominated solutions produced by the algorithm significantly dominate the ones generated by the current practice, in terms of fuel costs, while achieving similar path durations. The average number of non-dominated paths is 2.7, which allows decision-makers to ultimately visually inspect the proposed alternatives. 
\end{abstract}


\section{Introduction}\label{sec:intro}

Long-haul truck transportation is concerned with freight transportation from shipments' origins to destinations, with vehicle trips lasting from some hours to several days. Drivers performing long-haul transportation are subject to strict rules derived from Hours of Service (\hos{}) regulations. The aims of such regulations are to protect drivers and promote road safety by preventing accidents related to excessive fatigue. Therefore, \hos{} regulations typically limit the daily and weekly driving and duty times.

There exists a large body of literature integrating \hos{} regulations within long-haul transportation (see the literature review in Section \ref{sec:litHOS}).
The optimization problems in this context generally deal with routing and scheduling decisions aimed at determining where a driver should stop (for visiting customers or resting) and how long a rest should be.
Given the length of the routes in long-haul transportation, vehicles may need to refuel on several occasions. The overwhelming majority of the literature on long-haul transportation ignores refueling decisions and treats fuel costs as proportional to the traveled distance. Thus, implicitly assuming that fuel costs are uniform throughout the road network, and that refueling can be performed without any route deviations.
However, in practice fuel prices may differ considerably between countries (\cite{santos2017road}). For example, as derived from the \cite{eucom-oil} weekly oil bulletin, during 2019 the average diesel price in Italy was 17\% higher than in Germany. Moreover, fuel prices may be substantially different within the same urban area. For example, on the 12$^{\text {th}}$ of June 2020 the diesel prices in the area of Milan ranged between 1.148 and 1.809 \euro{} per liter (see \cite{MISE}). 

\begin{figure}[ht]
 \centering
 \includegraphics[scale=0.3]{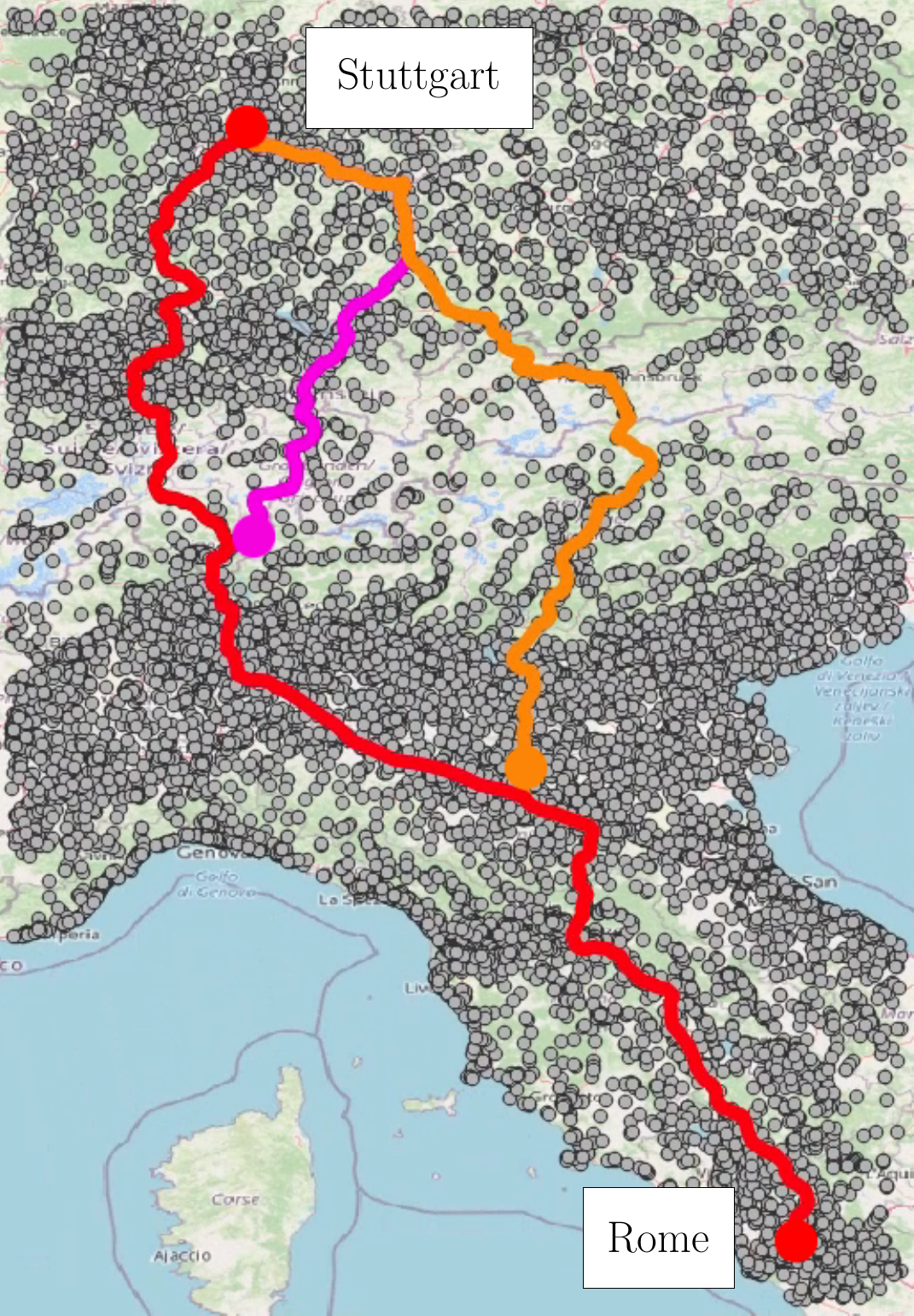}
 \caption{Refueling stations on paths from Rome to Stuttgart}
 \label{fig:pathRS}
\end{figure}

Despite the fact that fuel costs are a major cost component in transportation operations, the literature accounting for variable fuel prices in transportation planning decisions is very limited (\cite{neves2020multi}), 
as detailed in Section \ref{sec:litRef}. Most of the related papers have considered a limited number of refueling stations. For example, the heuristic proposed by \cite{bernhardt2017truck} was tested on graphs with up to 353 refueling stations. The motivation for considering a limited number of refueling stations in long-haul transportation stems from the assumption that vehicles would predominately refuel along highways. We argue that deviations from highways, to rural and possibly urban refueling stations, should be considered, as the additional distance may be offset by fuel cost savings. Therefore, the theoretical set of refueling stations to consider may be rather large, even in an urban area. For instance, in the metropolitan area of Milan (about 1575 km$^2$) there are 831 refueling stations (\cite{RELOMB}). A recent study by the European Commission shows that the average distance traveled at the EU level for road transportation is between 300 and 999 km (see \cite{eurostat}).
Therefore, the number of potential refueling stations is significantly amplified when considering a long-haul origin-destination trip. For example, the shortest path (based on travel time) of a trip from Rome to Stuttgart is about 1075 km. There are nearly 1090 refueling stations within a 5 km radius of this path (see Section \ref{sec:stopping} for the details on the refuel location search along a path). However, a deviation from the original planned path to a refueling station may lead the vehicle to proceed on a different path, after departing from a refueling station. An example of this is illustrated in Figure \ref{fig:pathRS}. Considering the origin-destination path from Rome to Stuttgart the shortest path is highlighted in red. A refueling stop at the orange refueling station may yield a modification in the original origin-destination path, such that the orange path is followed from the orange refueling station to the destination. This modification is due to the fact that the orange path is shorter than the path returning from the orange refueling station to the red path and proceeding to the destination. Such modifications may occur in many situations, e.g., the purple modification due to visiting the purple refueling station. Thus, the theoretical number of refueling stations can be extremely large. Indeed, there are 16,249 refueling stations in Figure \ref{fig:pathRS}. Fully accounting for such a number of nodes on a road network is impractical. This challenge is further exacerbated when considering the road network distance between a pair of nodes, and not their Euclidean distances. In fact, determining the distance between two points in a road network requires querying a geographic information system (\gis{}), and this might be cumbersome when the number of \gis{} calls increases, as well as expensive when the company does not own a \gis{} license.

In this paper we introduce the \name\ (\shortname). This problem originates from a collaboration with a company, Multiprotexion srl, specialized in the security of trucks that offers to its customers fleet planning and optimization services. The goal of the \shortname{} is to determine a path between an origin-destination pair complying with \hos{} regulations and fuel constraints. Two objectives are considered: the minimization of the fuel cost and the minimization of the duration of the path. While the first objective has a clear motivation in practice, the second is related to the optimization of the service level offered to customers, i.e., the faster is the transportation, the better is the service level. Thus, the resulting problem is a bi-objective optimization problem. The aim is to identify the Pareto frontier of feasible non-dominated paths. We note that the issue of finding alternative paths and routes in long-haul transportation applications is receiving increasing attention, e.g., \cite{caramia2009heuristic}.

We propose a heuristic solution algorithm for the \shortname{} that builds a set of non-dominated paths while interacting with a \gis{}. The algorithm dynamically determines a set of refueling locations and a set of rest locations that are of interest. Such locations are identified as being `on the way', i.e., they do not require a long detour from a considered path. The core idea of the algorithm is the following. A set of paths from the origin to the destination is devised and each path is individually explored starting from the origin. Once a stop is required, a set of feasible stop locations are considered based on the status of the vehicle and the driver. A path from each of those stop locations to the destination is computed. The process is repeated until the vehicle can reach the destination without violating \hos{} regulations and without running out of fuel. Results are presented on a set of instances where the origin-destination distance ranges from $500$ to $1500$ km and considering various initial \textit{stati} for the driver rest and the vehicle tank level.

The main contributions of the current work can be summarized as follows:

\begin{enumerate}
\item We introduce and define the \name.
\item We introduce a heuristic algorithm that interacts with a \gis{} to determine the relevant rest locations and refuel locations and to derive road network travel times between the considered nodes.
\item We define a set of experiments on instances based on road networks and show that the algorithm is capable of efficiently handling real-sized instances.
\item We compare the algorithm with a policy mimicking the current practice of a logistics company \item The results show that the algorithm is capable of determining a good variety of non-dominated solutions, on one side, and provides substantial fuel savings with respect to what is done in practice, on the other side. 
\item We highlight a number of managerial implications based on the obtained results.
\end{enumerate}

The remainder of the paper is organized as follows. In Section \ref{sec:lit} a brief review of the relevant literature is presented, in Section \ref{sec:problem_def} the \shortname{} is defined. Section \ref{sec:algorithm} describes the algorithm that is proposed for solving the \shortname{}. Computational experiments are presented in Section \ref{sec:exp}. Specifically, Section \ref{sec:impdet_instgen} reports the implementation details. The instance generation procedure is presented in Section \ref{sec:instances} and the 
 results are provided in Section \ref{sec:results}. Conclusions and managerial implications are drawn in Section \ref{sec:conclusions}.

\section{Literature review}\label{sec:lit}

In this section, the literature related to the \shortname{} is introduced. In Section \ref{sec:litHOS} we give a brief overview of the literature on optimizing long-haul transportation. This literature overwhelmingly ignores refueling decisions. In Section \ref{sec:litRef} we review the literature dealing with general routing problems and refueling decisions. Section \ref{sec:litRoad} is devoted to routing problems on road networks.

\subsection{Hours of service} \label{sec:litHOS}

The transportation science literature dealing with \hos{} regulations can be broadly categorized as i) long-haul vehicle scheduling, and ii) long-haul vehicle routing and scheduling. In the former category, a given sequence of (customer) locations should be visited by a vehicle. The problems in this category primarily schedule where the driver should stop and for how long a rest should be taken, in accordance with \hos{} regulations (e.g., \cite{AS2009}, \cite{Goel2009}). The fixed sequence of locations assumption is relaxed in the latter category, thus, simultaneously optimizing routing and scheduling decisions.

The Truck Driver Scheduling Problem (TDSP) was first addressed by \cite{Xu2003}, and then formally introduced by \cite{AS2009} (under the name of the \textit{trip scheduling problem}). Given a fixed sequence of locations to visit with time windows, \cite{AS2009} propose a polynomial algorithm that produces a feasible solution to the TDSP. A more efficient polynomial algorithm for similar problem settings is proposed by \cite{goel2012truck}. Both previously mentioned studies consider the TDSP subject to \hos{} regulations of the USA. Extensions to the rules associated with the legislation of other countries can be found in \cite{Goel2012JOSb,goel2012australia,Goel2012,goel2010europe}.

When the sequence of nodes to be visited is not established a-priori, we have the Vehicle Routing Truck Driver Scheduling Problem (VRTDSP). This problem has been extensively studied (e.g., \cite{Goel2009}, \cite{Kok2010}, \cite{Prescott2010}, \cite{Rancourt2013}, \cite{Goel2014}, \cite{Goel2017}). A comprehensive survey of VRTDSP solution methodologies is presented in \cite{tilk2020bidirectional}. \cite{kocc2018long} account for fuel costs in the VRTDSP in the context of idling options. Idling refers to the practice of leaving the vehicle's engine on during breaks to maintain a comfortable temperature or to use amenities such as television. However, the authors do not consider the decision about where refueling operations should take place in order to minimize fuel cost. In fact, in practical applications, fuel cost might differ remarkably on the basis of the location where refueling is done, especially in the case of international transportation. Furthermore, to the best of our knowledge, the VRTDSP literature does not account for distances on a road network.

\subsection{Refueling}\label{sec:litRef}

 The issue of heterogeneous fuel prices has been mainly addressed in a fixed route context. Given a fixed sequence of nodes to visit, the Fixed Route Vehicle Refueling Problem (FRVRP) is the problem of determining the sequence of refueling stops and the refueling amount for each stop in order to minimize the refueling cost. \cite{suzuki2014variable} proposes a pre-processing technique for the FRVRP to reduce the problem size without eliminating any optimal solution. A total of 16 instances with up to 495 refueling stations are solved with CPLEX. \cite{suzuki2014dss} consider the negative impact of carrying excessive amounts of fuel in the vehicle's tank. A formulation allowing the vehicle to retain some empty space in the tank is proposed for the FRVRP and tested on 24 instances with up to 199 refueling stations. \cite{lin2007linear} propose a linear-time algorithm for finding optimal vehicle refueling policies. The FRVRP is treated as a special case of the inventory-capacitated lot-sizing problem.

Refueling and routing decisions have also been studied jointly. \cite{suzuki2012decision} proposes a decision support system to tackle the traveling salesman problem with time windows and refueling. The author proposes a heuristic algorithm that sequentially solves a traveling salesman problem with time windows and a FRVRP. The results are reported on instances with up to 20 customers and with a density of refuel locations up to around one station every 39 kilometers on each arc. \cite{khuller2011fill} study the Variable-Route Vehicle-Refueling Problem (VRVRP) and propose a polynomial-time approximation algorithm for it. \cite{suzuki2013decision} study the bi-objective VRVRP with the aim of allowing carriers to jointly minimize fuel costs and vehicle mileage. \cite{bousonville2011vehicle} propose an extension of the Vehicle Routing Problem with Time Windows (VRPTW) to include refueling decisions. A fuel optimization model is developed and embedded in an insertion heuristic for the VRPTW. The proposed solution method is tested on instances with up to 100 customers and up to 441 refueling stations. \cite{neves2020multi} study a multi-period vehicle routing problem with refueling decisions where detours to reach refueling stations with lower prices are considered in order to minimize costs. A formulation is proposed for the problem and a branch-and-cut algorithm is presented for its solution. Computational results are presented on instances with up to 40 customers and 6 refueling stations.

To the best of our knowledge, only one report considers refueling together with \hos{} regulations. \cite{bernhardt2017truck} study the truck driver scheduling problem with rest periods, breaks and vehicle refueling in the context of international freight transportation. In particular, given a route and a set of refueling stations with different fuel prices along the route, the decisions relate to determining a time window to visit each customer, refueling stations to visit, refueling amounts and driver activities at stop locations. The latter include rest and refueling. The problem is a multi-criteria optimization problem with the goal of minimizing lateness, traveling time, and fuel expenditures. A Mixed Integer Linear Programming (MILP) model is proposed for the resulting problem, together with a pre-processing reduction technique to reduce the considered number of refueling stations. The proposed approach is tested on instances with up to 515 refueling stations. The main difference between the work of \cite{bernhardt2017truck} and the current paper is related to the representation of the underlying network. In fact, \cite{bernhardt2017truck} consider a complete graph where vertices represent customer and refueling station locations. Instead, we work directly on the road network and build the set of fuel and rest locations dynamically, as described in Section \ref{sec:algorithm}. 
\subsection{Road networks}\label{sec:litRoad}

The literature on routing problems defined on road networks is growing (\cite{GARAIX201062}, \cite{LETCHFORD2014331}, \cite{BenTicha2018}, \cite{BenTicha2019}, \cite{BENTICHA2019113}). The first issue to be considered is how to deal with the representation of the road networks. One possibility is to represent the network through a graph or a multi-graph where nodes and arcs represent elements of the road network (\cite{GARAIX201062}, \cite{LETCHFORD2014331}, \cite{BENTICHA2019113}). In particular, several multi-objective transportation planning problems based on GIS have been proposed (e.g., \cite{chen2008multi}). For a recent overview on optimization applications through GIS the reader is referred to \cite{murray2021contemporary}. This solution, however, might be not viable for long-haul transportation where the underlying network refers to a wide area, as the size of the corresponding graph (or multi-graph)
grows extremely fast with respect to the size of the represented area. In this case, it is preferable (and is often necessary) to deal directly with the road network by devising a procedure that interacts with a \gis{} and dynamically obtains the information needed through proper queries. This is the approach used in the current work.

\section{Problem definition}\label{sec:problem_def}

The \name\ (\shortname) is defined as follows. We are given an origin-destination pair, where the origin is denoted as $o$ and the destination as $d$. The vehicle has to transport goods from $o$ to $d$ through a path. A set of refueling stations $G$ is dynamically constructed, where each station $g \in G$ is associated with a price $p_g$ that is the unitary (per liter) fuel price. In addition, a set of rest locations $L$ is dynamically constructed, where a rest location corresponds to a place where the vehicle may stop and the driver may take a rest. In the following, we use the term \textit{stop location} to indicate either a refueling station or a rest location. The two considered objectives pertain to the refueling cost and the route duration. The former corresponds to the total cost of fuel consumed while traveling, whereas the latter consists of the driving time and rest periods. These are derived by assuming that the driver is subject to the European Union \hos{} regulations (EU: Regulation (EC) No 561/2006 and Directive 2002/15/EC, see \cite{eucom}).
In particular, we consider the following rules:
\begin{enumerate}
	\item \textit{Continuous driving rule}: After a maximum of four and a half hours of continuous driving, the driver has to take a break of at least 45 minutes. In the following, this is referred to as a \textit{break}.
	\item \textit{Maximum daily driving rule}: After a maximum of nine hours of driving, the driver has to take a rest of at least 11 hours. In the following, this is referred to as a \textit{daily rest}.
	\item \textit{Maximum weekly driving rule}: After a maximum of 56 hours of driving, the driver has to take a rest of at least 45 hours. In the following, this is referred to as a \textit{weekly rest}.
\end{enumerate}

These rules imply the necessity for the driver to stop to rest. Four stop types are therefore defined: for fuel, for a break, for a daily rest, and for a weekly rest. Note that the rules imposed by EU \hos{} regulations are more complex and include special cases. However, as a starting point, in this study we consider the subset of rules listed above. Nonetheless, our solution methodology can be adapted to more general \hos{} regulations. We note that, given the cost structure of our problem, in particular, the minimization of the tour duration, the following two additional rules are implicitly accounted for:
\begin{itemize}
	\item Within each period of 24 hours after the end of the previous daily rest period, a driver shall have taken a new daily rest period.
	\item A weekly rest period shall start no later than 144 hours after the end of the previous weekly rest period.
\end{itemize}

A traveling time $t_{ij}$ and a distance $c_{ij}$ are associated with each pair of locations $i,j \in \{o,d\}\cup G \cup L$. We assume that the vehicle consumes a constant amount of fuel per kilometer, which is denoted by $\phi$. In addition, the initial state of both vehicle and driver are known when departing from $o$. We denote as $f_o$ the fuel level in the vehicle tank when departing from $o$ and the capacity of its tank as $\tau$.
The number of hours remaining before the driver has to take a break, a daily rest, and a weekly rest, are denoted as $b_o$, $r_o$, and $w_o$, respectively.

A path from $o$ to $d$ is feasible if it satisfies the \hos{}
regulations and is such that the vehicle never runs out of fuel. The \shortname\ is an optimization problem where the path from $o$ to $d$ has to be determined together with the stops for rest and refueling. The two considered optimization objectives are:

\begin{enumerate}
	\item The minimization of the refueling cost. This corresponds to the total cost of fuel consumed while traveling from $o$ to $d$. If we denote by $\tilde{G}=\{g_0,g_1,...g_n\}$ the set of refueling stations visited along the path (with $g_0$ being the origin $o$ and $g_n$ being the destination $d$), then the refueling cost corresponds to $\sum_{i=0}^{n-1}\phi c_{g_ig_{i+1}}p_{g_i}$. Note that station $g_0$ is associated with a fuel price that corresponds to the average fuel price over all stations.
	
	\item The minimization of route duration. This is comprised of the driving time, the time required for breaks, daily rests, and weekly rests, as well as the refueling time.
	
\end{enumerate}

Given that long-haul transportation applications typically involve large distances, the theoretical size of $G$ and $L$ can be extremely large. Furthermore, considering a fully connected network, as is often the case in the mathematical programming formulation of transportation problems, yields an impractical problem size. The Rome to Stuttgart example presented in Section \ref{sec:intro} would imply computing a travel time matrix of about $264,030,001$ entries. An input of this size is unmanageable by state-of-the-art solvers, such as CPLEX. Indeed, some preliminary experiments showed that it was not even possible to load the corresponding problem formulation in CPLEX, on a 64GB workstation. Also, this might be expensive in practice for companies that do not own a routing software and have to pay for querying a \gis{} software. Finally, computing the entire set of paths is computationally cumbersome. Therefore, we adopt a heuristic solution approach, which is presented in the subsequent section.
The heuristic reduces the set of rest and refuel locations considered to those that are more likely to be part of non-dominated routes.

\section{Solution algorithm on a road network}\label{sec:algorithm}

The algorithm works by iteratively building a set of feasible $o-d$ paths. Considering this set, Pareto optimal paths are found.
We define a \textit{temporary path} as a path between a location $\bar{o} \in \{o\}\cup G\cup L$ and $d$, which may be feasible or not. Initially, the algorithm constructs a set of temporary paths between $o$ and $d$. For each infeasible path, stop locations are added and combined, and the path is updated in order to generate a feasible sub-path between $o$ and each inserted location. Then, temporary paths are created between each inserted location and $d$, and the process is repeated until a feasible path is obtained. 

When going from $o$ to $d$, the vehicle performs two activities: it travels or it stops. Because of the huge size of the set of stop locations on a road network, the sets $G$ and $L$ are heuristically built by progressively exploring the generated paths and querying the \gis{} for the required stops. We define a \textit{stop location} as the position where the stop takes place and we identify four \textit{stop types}:

\begin{itemize}
	\item $F$: stop for refueling;
	\item $B$: stop for a break;
	\item $D$: stop for a daily rest;
	\item $W$: stop for a weekly rest.
\end{itemize}

A daily rest (stop of type $D$) is also a break (stop of type $B$). This means that, after a stop of type $D$, the driver has a maximum availability in terms of daily driving time of nine hours, and continuous driving time of four and a half hours. Similarly, a weekly rest (stop of type $W$) is also a daily rest (stop of type $D$) and, therefore, also a break (stop of type $B$). A stop of type $F$ can be combined with a stop of type $B, D$ or $W$. When this happens, the refueling is assumed to take place during the rest. Otherwise, a fuel stop is assumed to last 15 minutes. Refueling operations are assumed to always fill the tank. This is consistent with what is observed by the company which inspired this work and with carriers' practice.
A \textit{stop} is defined by its location and the subset of stop types that define the operations to be performed at the stop, e.g., taking a daily rest and refueling.

\begin{figure}[ht]
	\centering
	\includegraphics[width=\textwidth]{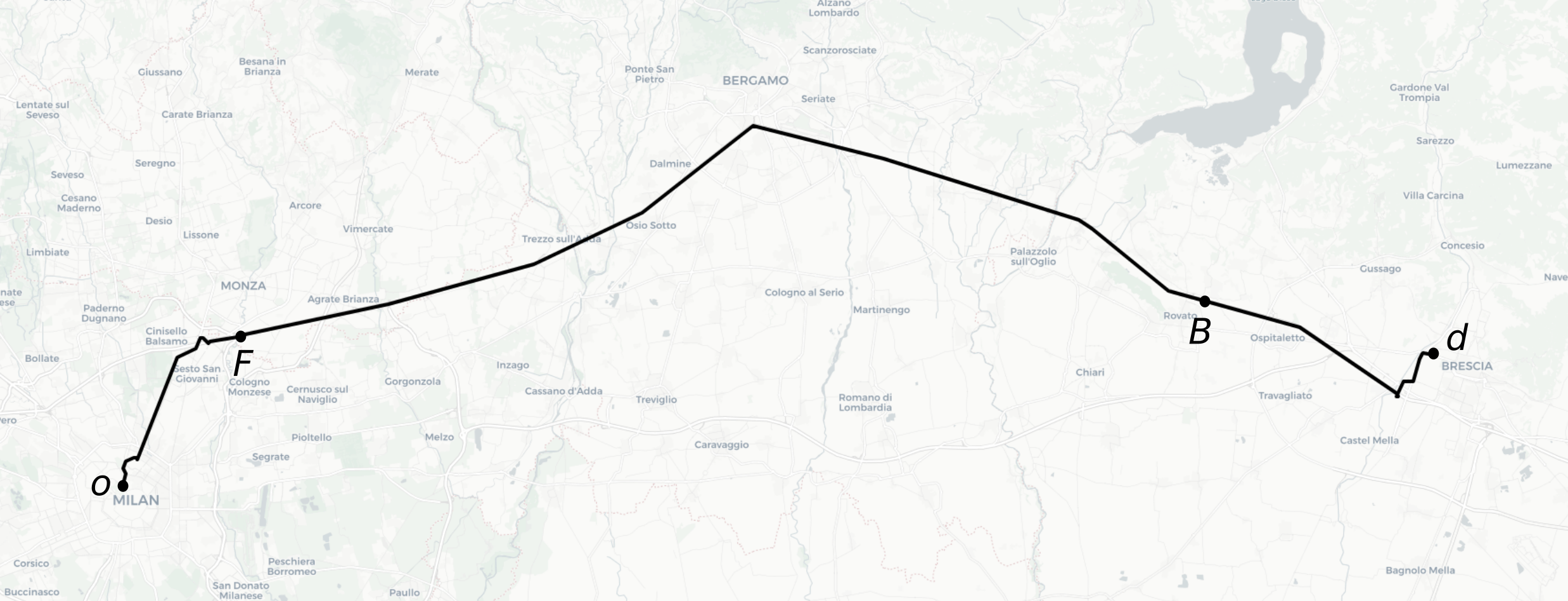}
	\caption{A path from Milan to Brescia}
	\label{fig:path}
\end{figure}

A path between $o$ and $d$ may include multiple stops. We call an \textit{arc} a portion of the path between two nodes of $\{o,d\}\cup G \cup L$. Thus, a path is a concatenation of consecutive arcs. Each arc corresponds to a geographical trajectory retrieved by querying the \gis{}. A path between Milan and Brescia with two stops (one for refueling and one for a break), and hence three arcs, is shown in Figure \ref{fig:path}.

The solution algorithm for the \shortname\ is outlined in Table \ref{algo:treeExp}.
The first step of the algorithm (line 3) is to determine a set of $k$ temporary paths from $o$ to $d$. These correspond to the $k$ fastest paths (in terms of travel times), ignoring any \hos{} and refueling constraints. Such paths can be retrieved from standard pathfinding libraries such as \cite{gh}. Specifically, function \textit{branchToDestination} takes as input a node and provides as output the $k$ fastest (time) temporary paths from the node to destination $d$. We note that such temporary paths from $o$ to $d$ are likely to be infeasible in long-haul transportation.

\begin{table}[h]
\begin{algorithm}[H]
	\small
	\caption{Solution algorithm} 
	\begin{algorithmic}[1]
		\STATE \textbf{Input:} $o,d,k$
				\STATE \textbf{branchToDestination(}$o$, $k$\textbf{)} \label{line:beginning}
		\STATE infeasibleArc := \textbf{selectInfeasibleArcToDestination(}tree\textbf{)} \label{line:siatd}
		\WHILE{ infeasibleArc != null}
		\label{line:loop}
		\STATE $S=\emptyset$
		\STATE $\bar{o}$ := origin of infeasibleArc
		\STATE $Q$ := \textbf{findStopLocations($\bar{o},d$)} (Section \ref{sec:stopping})
		\FORALL{$q \in Q$}
		\STATE $S$ := $S$ $\cup$ \textbf{combineStopTypes(}$q$\textbf{)} (Section \ref{sec:routes})
		\ENDFOR
		\FORALL{$s \in S$}
		\STATE create arc from the $\bar{o}$ to $s$
		\STATE \textbf{branchToDestination(}$s,k$\textbf{)}
		\ENDFOR
		\STATE remove infeasibleArc $(\bar{o},d)$
		\STATE infeasibleArc := \textbf{selectInfeasibleArcToDestination(}tree\textbf{)}
		\ENDWHILE
		\STATE Identify Pareto optimal paths \label{line:Pareto}
		\STATE \textbf{Output:} Pareto optimal paths
	\end{algorithmic}
\end{algorithm}
\caption{The pseudo-code of the solution algorithm} \label{algo:treeExp}
\end{table}

After constructing the initial $k$ temporary paths, the algorithm treats each one independently. For each infeasible temporary path, a set of stops are added with the aim of generating feasible paths. During this process, several other temporary paths are generated from the newly identified stops to the destination.
The overall construction of paths follows a tree-based procedure where the tree is rooted in $o$, each intermediate node is associated with a location in $G \cup L$ and the operation performed in it, and all leaves correspond to $d$. Each arc of the tree represents an arc of a path. When the construction phase is finished, each path on the tree from the root ($o$) to a leaf ($d$)
corresponds to a feasible path.

The tree is initialized in line \ref{line:beginning} by setting $o$ as the root of the tree, with $k$ branches corresponding to the $k$ temporary $o$ - $d$ paths and with $k$ leaves all corresponding to $d$. The function \textit{selectInfeasibleArcToDestination} in line \ref{line:siatd} and the loop in line \ref{line:loop} looks for an infeasible arc in the tree. If no infeasible arc is found, the algorithm stops as all generated paths are feasible.
	
Otherwise, an infeasible arc of the tree is selected according to a depth-first search algorithm. Considering an infeasible arc $(\bar{o},d)$, stop locations are searched along its associated temporary path by function \textit{findStopLocations} (explained in Section \ref{sec:stopping}). Considering each of these stop locations, the function \textit{combineStopTypes} creates stops (specified in terms of location, type and duration) according to the rules explained in Section \ref{sec:routes}. This procedure entails duplicating stop locations in order to allow different combinations of fueling and resting. Each found stop $s$ is inserted in the tree, and the arc ($\bar{o}$,$s$) is added to the tree. The travel time of this arc corresponds to the fastest path between $\bar{o}$ and $s$. The infeasible arc $(\bar{o},d)$ is removed from the tree. Then $k$ temporary paths between $s$ and $d$ are created through the \textit{BranchToDestination} function, and appended to the tree.

\begin{figure}[ht]
	\centering
	\includegraphics[scale=0.6]{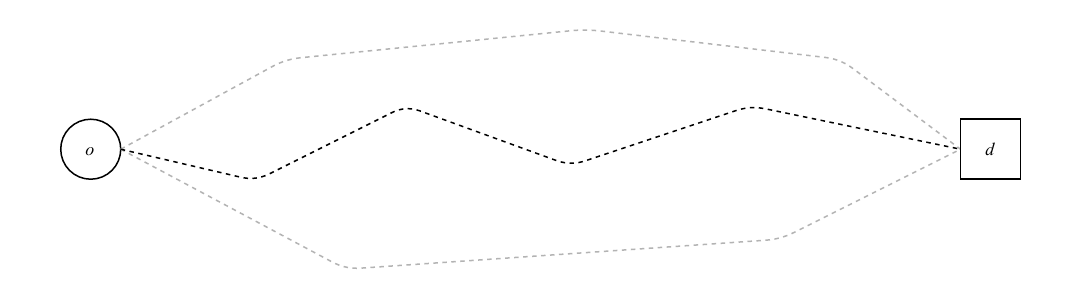}
	\includegraphics[scale=0.6]{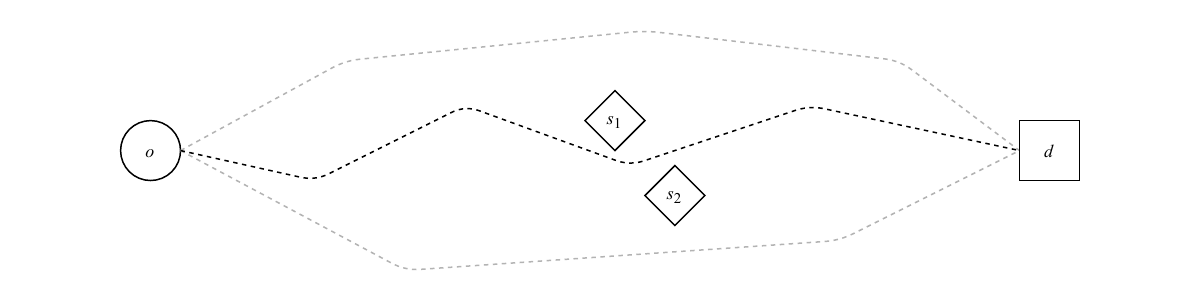}
	\includegraphics[scale=0.6]{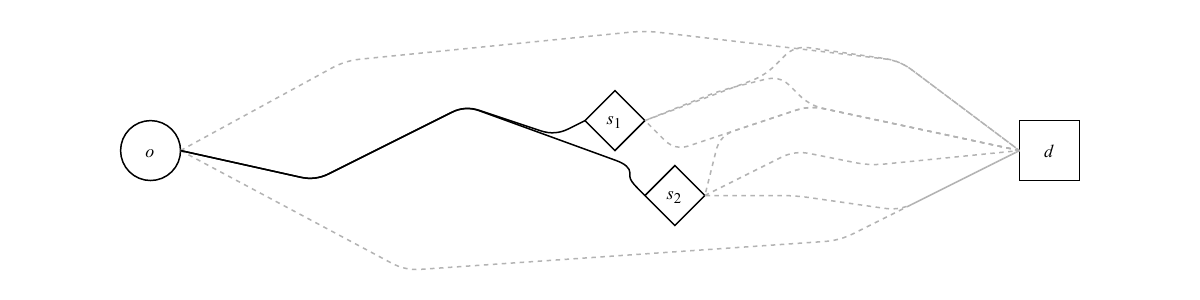}
	\caption{Tree generation example with $k=3$}
	\label{fig:pathtree}
\end{figure}

A simplified example of the previously described procedure, without the combine stop types option, for $k=3$, is provided in Figure \ref{fig:pathtree}. In the top panel, three temporary paths from $o$ to $d$ are identified as dashed lines. The selected infeasible arc is black while the others are gray. In the central panel, stop locations are found along the path. In the bottom panel, the feasible arcs from the origin to the stop locations are marked as solid black lines and the temporary paths from each stop to the destination $d$ are marked as grey dashed lines. Note that, even if paths may overlap (as in the bottom panel of Figure \ref{fig:pathtree}), they are associated with different arcs in the tree.

Each node of the tree is a stop. Thus, it denotes the stop location and stop type, from which we derive the resulting fuel level and driver resting conditions. Furthermore, all arcs of the tree are feasible with the only possible exceptions being the ones ending at destination $d$. Thus, when exploring the tree searching for an infeasible arc, only arcs going to $d$ are considered, that is, only arcs ending in a leaf of the tree. Finally, in line \ref{line:Pareto}, a Pareto frontier is created, by examining all generated feasible paths based on the two objective functions: fuel cost and route duration.

In the next sections, we first describe how the sets $G$ and $L$ are generated (Section \ref{sec:stopping}), then describe how stops types are combined and feasible routes are constructed (Section \ref{sec:routes}).
Table \ref{tab:notation} summarizes the notation used in the following sections for the description of the algorithmic procedures.

\begin{table}[ht]
	\centering
	\begin{tabular}{cl}
		Notation &\multicolumn{1}{c}{Description} \\
		\hline
		\hline
		\multicolumn{2}{c}{Stops}\\
		\hline
		$s_F$& Stop for refueling\\
		$s_B$& Stop for making a break\\
		$s_F$& Stop for making a daily rest\\
		$s_W$& Stop for making a weekly rest\\
		$s_{FB}$& Stop combining refueling and break\\
		$s_{FD}$& Stop combining refueling and daily rest\\
		$s_{FW}$& Stop combining refueling and weakly rest\\
		\hline
		\multicolumn{2}{c}{Radius search}\\
		\hline
		$C$& Center of the search\\
		$\rho$& Radius of the search\\
		$I$& Length of the interval of the search\\
		\hline
		\multicolumn{2}{c}{Refuel location search}\\
		\hline
		$\beta\tau$& Upper bound on tank level for refuel location search\\
		$\gamma\tau$& Lower bound on tank level for refuel location search\\
		$C_{\beta\tau}$& Location reached when the tank level is $\beta\tau$\\
		$C_{\gamma\tau}$& Location reached when the tank level is $\gamma\tau$\\
		\hline
		\multicolumn{2}{c}{Rest location search}\\
		\hline
		$\alpha$& Next point in time when a stop due to \hos{} regulations is needed\\
		$\delta$& Threshold defining the beginning of the rest location search\\
		$\alpha\delta$& Point in time when the rest location search starts\\
		$C_{\alpha\delta}$& Location reached at time $\alpha\delta$\\
	\end{tabular}
	\caption{Summary of notation}
	\label{tab:notation}
\end{table}

\subsection{Defining the set of refuel and rest 	locations}\label{sec:stopping}

In this section, the \textit{findStopLocations} function, which finds stop locations on the path corresponding to an infeasible arc $(\bar{o},d)$, is described. 
Regarding the search for refuel locations, two points corresponding to fuel level $\beta\tau$ and $\gamma\tau$ are found along the path and refuel locations are looked for in the interval defined in such a way.
Regarding the search for rest locations, we denote by $\alpha$ the next point in time when the driver has to stop due to \hos{} regulations, whether it is due to a break, a daily stop, or a weekly stop being required. The search for a rest location is carried out starting from the point where the driver reaches $\delta \alpha$ driving time along the path. Choosing a value of the percentage $\delta$ not too close to $1$, e.g., $0.95$, allow us not to exceed driving time $\alpha$ in the search for a rest location on the road network. More details on the refuel location search and the rest location search are reported in the remaining part of this section.

Both searches are defined based on a \textit{radius search}, which has three inputs:

\begin{itemize}
	\item Center of the search ($C$) - \textit{location on the path}.
	\item Radius of the search ($\rho$) - \textit{maximum Euclidean distance}.
	\item Stop type - \textit{type of stop to be searched}.
\end{itemize}

\begin{center}
	\begin{figure}[ht]
		\centering
		\includegraphics[scale=0.3]{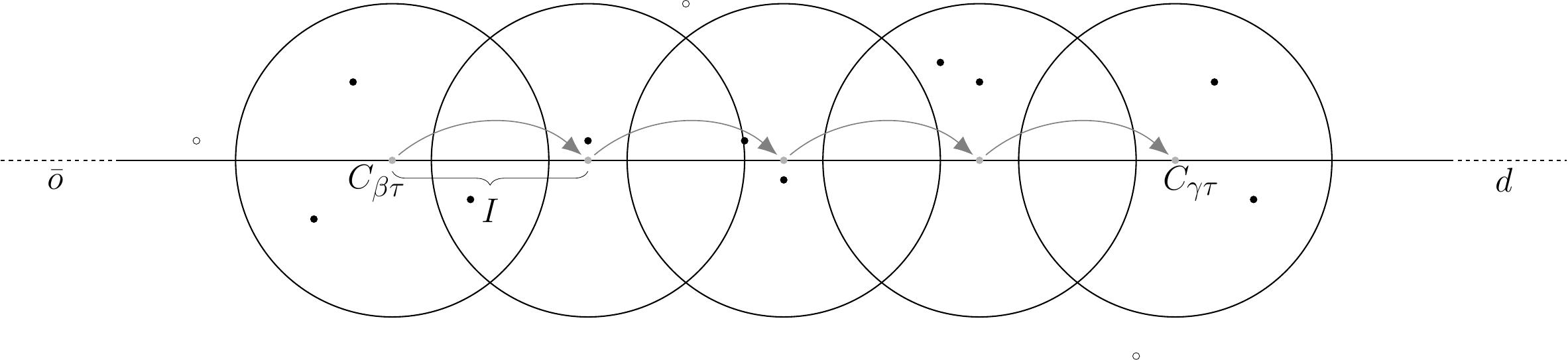}
		\caption{Refuel location search}
		\label{fig:intervalsearch}
	\end{figure}
\end{center}

\begin{figure}[ht]
	\centering
	\includegraphics[scale=0.3]{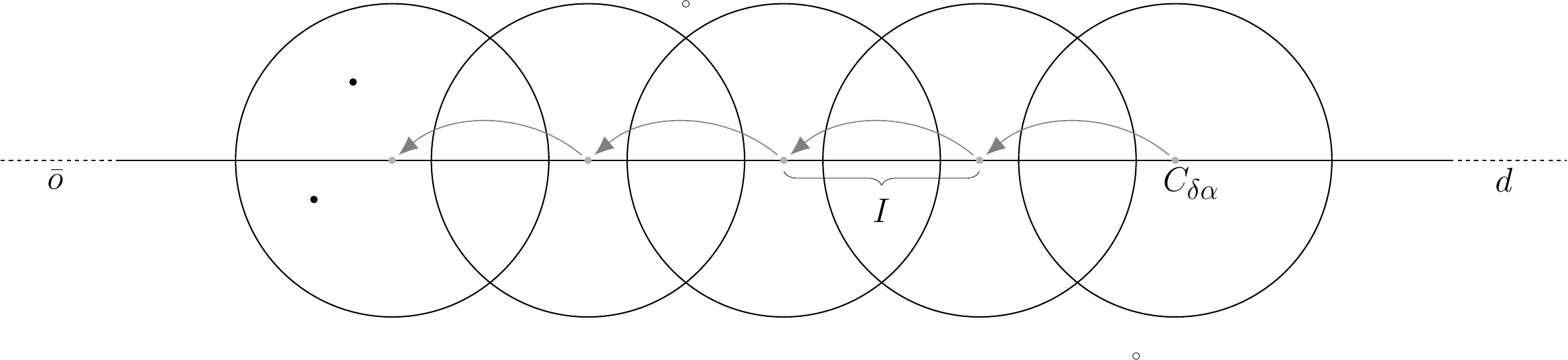}
	\caption{Rest location search}
	\label{fig:spotsearch}
\end{figure}

Given these inputs, the search for a feasible stop location is performed in the resulting area (see Figure \ref{fig:elementary}). The radius search is carried out differently according to whether a refuel location or a rest location has to be found. In particular:

\begin{itemize}
	\item \textit{Refuel location search}: when a refuel location must be found along a temporary path between $\bar{o}$ and $d$, the search is carried out between two locations $C_{\beta\tau}$ and $C_{\gamma\tau}$, defined on the basis of when the upper and lower bounds on the tank level are reached. The radius search is repeated on the path at intervals of length $I$ km (see Figure \ref{fig:intervalsearch}), where $I$ is a given parameter.
	\item \textit{Rest location search}: given a temporary path between $\bar{o}$ and $d$ and a location $C_{\delta\alpha}$ on that path (defined on the basis of when the 	specified amount of time left is reached), a radius search is performed on the path going backward at a specified interval of length $I$ km until at least one feasible rest location is found (see Figure \ref{fig:spotsearch}).
\end{itemize}

 The definition of the lower bound $C_{\beta\tau}$ (for the interval in which refuel locations are searched for) helps in restricting the number of potential stops considered. Indeed, if no lower bound is considered, then any refuel location along the path (or within a certain distance from the path) could be considered as a potential location and should be evaluated, thus making the number of evaluations excessively large.
 
Among the rest locations found, we keep the one with the cheapest insertion cost between $\bar{o}$ and $d$. The insertion cost is evaluated in terms of traveling time, i.e., we keep the path that leads to the smallest increase in traveling time. This choice is due to the fact that, when searching for a rest location, no refueling is made so the fuel price has no impact. Also, we assume that the distance traveled is proportional to time, so that the path associated with the smallest increase in traveling time is the one associated with the smallest increase in the distance traveled. A similar rule is used for reducing the number of refuel locations, based on the two objectives of the optimization. More precisely, a refuel location is not considered if there exists another refuel location with the same or cheaper fuel price and with a cheaper insertion cost. This rule limits the size of the search tree.

\begin{center}
	\begin{minipage}{0.48\linewidth}
		\centering
		\includegraphics[scale=0.55]{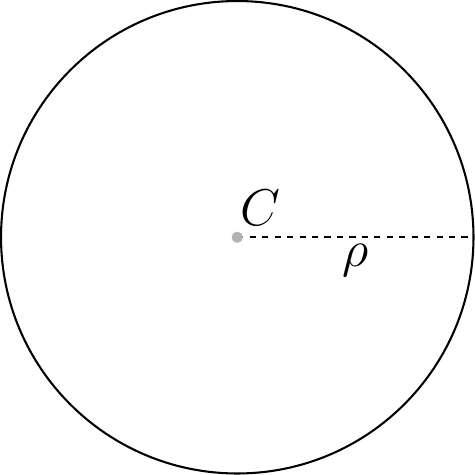}
		\captionof*{figure}{Search}
	\end{minipage}	\hfill
	\begin{minipage}{0.48\linewidth}
		\centering
		\includegraphics[scale=0.55]{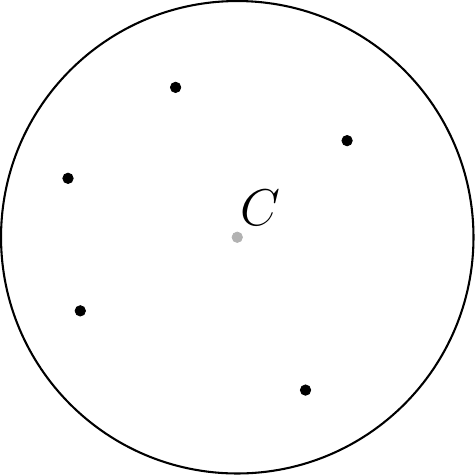}
		\captionof*{figure}{Result}
	\end{minipage}
	\captionof{figure}{Radius search} \label{fig:elementary}
\end{center}

\begin{figure}[ht]
	\centering
	\includegraphics[width=\textwidth]{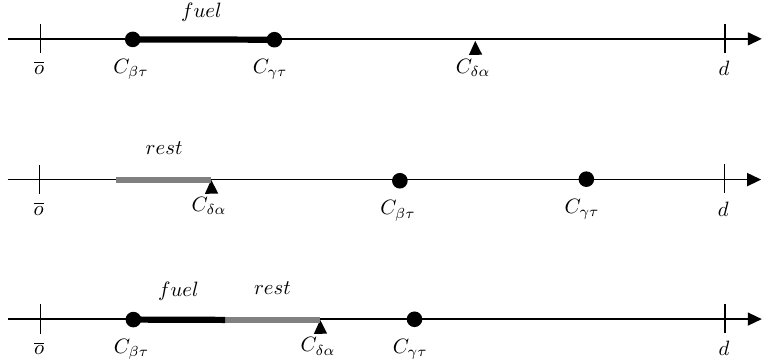}
	\caption{Stop location search 	}
	\label{fig:stopstragegy}
\end{figure}

Depending on the status of the fuel level and driving time, the execution of the \textit{findStopLocations} function may yield one of the following outcomes: i) refueling stop locations (when the first cause of infeasibility is due to fuel), ii) a rest stop location (when the first cause of infeasibility is due to \hos{} regulations), or iii) refueling stop locations and a rest stop location (when the arc is infeasible both with respect to fuel autonomy and \hos{} regulations, and the two searches overlap). These situations are depicted in Figure \ref{fig:stopstragegy}. Considering a path from $\bar{o}$ to $d$, the part of the path in which a fuel search is performed (the path between locations $C_{\beta\tau}$ and $C_{\gamma\tau}$ in Figure \ref{fig:intervalsearch}) is represented by a thick black line. The starting point of the backward rest location search (location $C_{\delta\alpha}$ in Figure $\ref{fig:spotsearch}$) is represented by a triangle, and the interval in which a rest location has been searched is represented by a thick gray line. In the case represented on the top panel, refueling is the first cause of infeasibility, as it is closer to $\bar{o}$. In the middle panel, rest is required before refueling is. In the bottom case, rest is required within the refueling interval. In this case, refueling stations are searched until the rest is requested. Then, once a rest is required, a combination of fuel and rest locations are searched. 

\subsection{Combining stops}\label{sec:routes}

 In this section we describe how the identified stop locations, from the previous section, are duplicated and converted into stops (specified in terms of location, type, and duration). We recall that, once a stop is determined, the resulting fuel level and the residual driving time, according to the \hos{} regulations, upon exiting the stop can be computed.

The algorithm has the capability to combine fuel and rest stops. Furthermore, it considers combining rest stops. For instance, if a break has to be taken because the driver used all the driving time, it might be advantageous to take a daily rest if the daily driving time remaining is almost over. This is particularly relevant when minimizing the route duration since combining such stops does not influence the fuel cost.

The function \textit{combineStopTypes}, used in the algorithm (see Table \ref{algo:treeExp}), considers two main possibilities of combining stop types for a given location: i) combining fuel with rest stops and ii) combining different rest stops. In principle, the options of combining stops are performed by duplicating the relevant stop locations.

Combining fuel with rest stops is done as follows. A fuel stop location $q\in G$ is replicated into four stops: refueling (only) and refueling with the three resting options. In particular, the following nodes are created in the tree:

\begin{enumerate}
	\item $s_F$: a stop of type $F$;
		\item $s_{FB}$: a stop where types $F$ and $B$ are combined;
		\item $s_{FD}$: a stop where types $F$ and $D$ are combined;
			\item $s_{FW}$: a stop where types $F$ and $W$ are combined.
		\end{enumerate}

The possibility of combining different types of stops when a rest is needed is performed as follows. When a stop location in $q \in L$ is identified, i.e., when a rest is required, stops $B$, $D$, and $W$ are evaluated. In particular, when a $B$ stop is required, stops $D$ and $W$ are also considered and when a $D$ stop is required, stop $W$ is also considered. Based on these rules, the following nodes are created in the tree:

\begin{enumerate}
	\item $s_B$: a stop of type $B$;
	\item $s_D$: a stop of type $D$;
	\item $s_W$: a stop of type $W$.
\end{enumerate}

\begin{figure}[ht]
	\centering
	\includegraphics[scale=0.8]{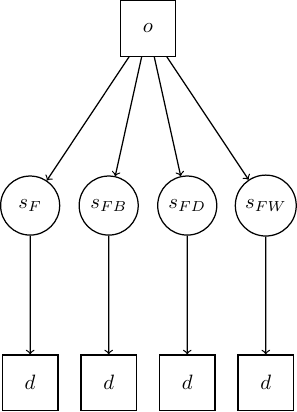}
	\caption{Tree-to-destination with fuel stops}
	\label{fig:fueltree}
\end{figure}

\begin{figure}[ht]
	\centering
	\includegraphics[scale=0.8]{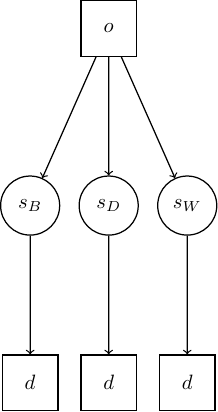}
	\caption{Tree-to-destination with rest stops when a break is required}
	\label{fig:stoptree}
\end{figure}

An example of the resulting tree is represented in Figure \ref{fig:fueltree} for a stop in a node in $G$ and in Figure \ref{fig:stoptree} for a stop in a node in $L$.

\section{Computational experiments}\label{sec:exp}

In this section we present the computational campaign. Specifically, Section \ref{sec:impdet_instgen} describes the details of the implementation of the algorithm illustrated in Section \ref{sec:algorithm}. The procedure for instance generation is presented in Section \ref{sec:instances} while Section \ref{sec:results} shows the results. The experiments are aimed at gaining insights regarding the problem difficulty and solution structure. In addition, a comparison with the current practice of the company is performed.

\subsection{Implementation details and instance generation} \label{sec:impdet_instgen}

The algorithm has been implemented in Java, while the base maps have been obtained through \cite{osm}. The considered paths between locations have been generated with \cite{gh}. Locations have been detected with \cite{op} (a web-based data filtering tool for \cite{osm}) by searching for nodes with the tag \textit{amenity} equal to \textit{fuel}, i.e., refuel locations. We assume that any type of stop can take place at a refuel location and that any refuel location can serve trucks. Furthermore, because of the lack of information on the availability of parking for trucks in accommodation facilities, we do not include any rest location in addition to refueling stations. In practice, this means that $L = G$. While this might not reflect reality, this assumption is consistent with the information available in the tools mentioned above. Note that the restricted availability of data is due to the fact that the current implementation is a prototype of the tool that the company motivating this work aims at developing. Being in its early stage of development, the first goal of the company was to verify whether such a tool might be beneficial for its customers (and the following results show that it is indeed beneficial). Thus, the company has not yet acquired a licensed map software with more detailed information. This forced us to turn to OpenStreetMap, which is free and thus ideal for prototype testing (despite lacking complete information). Moreover, from a computational perspective, this assumption entails that the algorithm is validated (as will be shown in Section \ref{sec:results}) on relatively large $G$ and $L$, and thus should be effective on sets with a smaller cardinality.
Fuel prices, at a country level, have been obtained from the Weekly Oil Bulletin of the European Commission (see \cite{eucom-oil}). The station fuel price was generated by perturbing the average national price up to $\pm15\%$, with the perturbation sampled uniformly from the resulting interval. This choice was made due to the absence of a consolidated accessible source of data on European daily fuel prices at fuel stations.

The number of temporary paths $k$ generated between $\bar{o}$ and $d$ was set to three. Out of $k$ we only consider paths which are less than $\sigma = 1.1$ times the shortest path between $\bar{o}$ and $d$.
With respect to the location search, fuel stops are searched in the surroundings of the path between locations $C_{\beta\tau}$ and $C_{\gamma\tau}$, where $\beta = 0.10$ and $\gamma = 0.05$, that is, between the point where the vehicle reaches 10\% of the fuel tank and the point where this measure reaches the 5\%. Rest stops are searched from the point $C_{\delta\alpha}$. The value was set to $\delta = 0.05$, that is when the driver has 5\% of the allowed driving time left since the last stop. To avoid extremely short time intervals, the value $\delta\alpha$ was imposed to be 10 minutes minimum. The radius $\rho$ of the search and the interval $I$ of the search are both set to 5 km. These values apply for both the rest and the refuel location search. 

The travel time of each path was determined using the open-source version of \cite{gh}, which provides travel times related to car paths. We have adjusted speeds to 100 km/h to reflect the fact that trucks adhere to lower speed limits, compared to private vehicles on a highway. We note that this speed was assumed on all traveled distances. This assumption is necessary due to the difficulty of computing the travel times for various road segments in \cite{gh}. Furthermore, this assumption is reasonable given that, in long-haul transportation, trucks predominantly use highways. The distances are obtained by multiplying travel times by the speed. We note that using navigation software with heavy-duty vehicle specifications may yield more accurate travel times. This, however, will not influence the implementation of the algorithm. 

In practice, one may obtain precise truck paths through platforms like \cite{gh} by paying a cost depending on different factors. The resulting paths, in this case, would consider roads that explicitly permit access to a given truck type, while accounting for its specific speed limit.
The cost of such queries depends on several parameters and is difficult to assess. We have therefore reported the number of calls to \cite{gh} in Section \ref{sec:results} as an indicator of the potential \textit{cost} of queries to a user not owning a truck navigation software.

\subsection{Instance generation}\label{sec:instances}

Instances have been generated as follows. Three values have been chosen for the initial temporary path length, namely, 500, 1000, and 1500 km. Three $o-d$ pairs have therefore been selected such that the length of the fastest path would be close to the corresponding path lengths.
These pairs are:

\begin{itemize}
	\item $500$ km: Freiburg im Breisgau (DE) to Maastricht (NL);
	\item $1000$ km: Paris (FR) to Brescia (IT);
	\item $1500$ km: Montreux (CH) to Timisoara (RO);
	\end{itemize}

The fastest and shortest paths length for the three $o-d$ pairs are reported in Table \ref{tab:fast_short}.
The shortest paths could include the so-called \textit{tertiary} roads, defined in OpenStreetMap as those roads with ``low to moderate traffic which link smaller settlements such as villages or hamlets'', which are not desirable when planning truck paths. Thus, we use the fastest paths in the proposed algorithm.

\begin{table}[ht]
	\centering
	\begin{tabular}{crr}
		Instance & \begin{tabular}[r]{@{}c@{}}Fastest\\ path (km) \end{tabular} & \begin{tabular}[r]{@{}c@{}}Shortest\\ path (km)\end{tabular} \\
		\hline
		Freiburg - Maastricht & 515.8 & 430.7 \\
		Paris - Brescia & 929.3 & 862.5 \\
		Montreux - Timisoara & 1488.1 & 1375.4 \\
				\hline
	\end{tabular}
	\caption{Fastest and shortest path (in km) for each $o-d$ pair}
	\label{tab:fast_short}
\end{table}

Three scenarios have been generated for the initial status of the driver:

\begin{itemize}
	\item b: just after a break, with half of the day and of the week hours remaining, i.e., $b_o=4.5$, $r_o=4.5$ and $w_o=28$;
	\item d: just after a daily rest, with half of the week hours remaining, i.e., $b_o=9$, $r_o=4.5$ and $w_o=28$;
	\item w: one work day remaining before a week rest, i.e., $b_o=4.5$, $r_o=9$ and $w_o=9$.
\end{itemize}

The capacity of the vehicle tank $\tau$ has been assumed to be 500 liters. The fuel consumption is assumed to be 3.5 km/l. This was established by averaging the fuel consumption for vehicles above 11.5t, as reported by the Italian Ministry of Transportation (see \cite{mit}). This value reflects the average characteristics of a long-haul truck in different load and road conditions (empty/full trailer, slope, wind, etc). Five scenarios have been tested with respect to the initial fuel level $f_o$: $\{10, 25, 50, 75, 100\}$ percent of the tank capacity. Thus, a total of 45 instances have been generated. Instances are referred to by reporting each characteristic separated by an underscore, e.g., the instance of 500 km where the driver has taken a break before departing and the vehicle has 10\% of the tank capacity at the origin is denoted as ``500\_b\_10''. 

\subsection{Computational results}\label{sec:results}
In this section the computational results for the generated instances are reported. Table \ref{tab:summary} provides a summary of the total number of feasible paths generated and the number of non-dominated paths found. Additional statistics on each instance are reported in the rightmost columns of the table, i.e., the run time of the algorithm, the percentage of run time ascribed to \cite{gh} (GH), \cite{op} (OP), and the algorithm, the number of calls to \cite{gh} (for one or multiple paths), 
the total number of non-dominated fuel and rest stop locations found. A first observation is that the number of non-dominated paths is very small with respect to the number of paths generated.
The number of non-dominated paths and the run time of the algorithm increase with the distance and tend to be lower when starting with a full or almost full tank.
A second observation is that a large percentage of the run times is ascribed to the handling of geographical information (i.e., \cite{gh} and \cite{op}). This result further highlights the cost to be paid when dealing with complex data such as that provided by \gis{}.

Tables \ref{tab:results}--\ref{tab:results1500_3} report the details on the non-dominated paths for each instance. In particular, Table \ref{tab:results} lists the results for the 500 km instances, Table \ref{tab:results1000} lists the results for the 1000 km instances. Table \ref{tab:results1500_2} lists the results for the 1500 km instances with b and d and Table \ref{tab:results1500_3} lists those for the 1500 km instances with w. We report, for each path, the length (kilometers), the fuel cost, the liters of fuel used, the total time, the travel time, the stop time (note that the total time corresponds to the sum of travel time and stop time), and the number of stops made for each type.

Furthermore, the performance of the algorithm is compared against an algorithm representing the current practice of the drivers (denoted as ``CP'' in the ``Path \#'' column). In the algorithm representing the current practice, drivers are assumed to behave myopically. In particular:
\begin{itemize}
\item Drivers are not assumed to consider alternative paths to the fastest one to the destination, i.e., $k = 1$.
\item Drivers look for a fuel location only shortly before refueling is required, that is, they consider refueling in locations between the points where 2\% and 1\% of the fuel capacity is reached. If two or more locations are found, the one with the cheapest distance insertion is considered. Rest locations are sought in the same way as the proposed algorithm.
\item Drivers are assumed to be able to combine refueling and resting only as follows: whenever the driver stops for fuel, if a \hos{} rest is required in 30 minutes or less, the fuel and rest stops are combined at the fuel stop location.
\end{itemize}

This representation of the current practice is justified by the knowledge gathered by the company on the behavior of the drivers.

\begin{table}[ht]
    \centering
\resizebox*{!}{1.4\columnwidth}{%
\begin{tabular}{cccrrrrrrrrr}
\multicolumn{3}{c}{Instance} &
\multicolumn{2}{c}{Paths}&
\multicolumn{6}{c}{Statistics}\\
\cmidrule(lr){1-3} \cmidrule(lr){4-5} \cmidrule(lr){6-12}
km & HoS & fuel & Total & \begin{tabular}[r]{@{}c@{}}Non-\\ dominated \end{tabular} 
& \begin{tabular}[r]{@{}c@{}}Total \\Run \\ Time(s) \end{tabular}
& \begin{tabular}[r]{@{}c@{}}\% \\ time\\ GH \end{tabular}
& \begin{tabular}[r]{@{}c@{}}\% \\ time\\ OP \end{tabular}
& \begin{tabular}[r]{@{}c@{}}\% time\\ Algorithm \end{tabular}
& \begin{tabular}[r]{@{}c@{}}\# calls \\ to GH \end{tabular}
& \begin{tabular}[r]{@{}c@{}}Evaluated \\ fuel stops \end{tabular}
& \begin{tabular}[r]{@{}c@{}}Evaluated \\ rest stops \end{tabular}
\\

\hline

500 & b & 10 & 296 & 2 & 245 & 68 & 28 & 4 & 3409 & 220 & 986 \\
500 & b & 25 & 132 & 5 & 89 & 47 & 50 & 3 & 1065 & 91 & 216 \\
500 & b & 50 & 6 & 1 & 5 & 78 & 20 & 2 & 79 & 0 & 29 \\
500 & b & 75 & 6 & 1 & 5 & 77 & 21 & 2 & 79 & 0 & 29 \\
500 & b & 100 & 6 & 1 & 5 & 76 & 21 & 3 & 79 & 0 & 29 \\
\cmidrule(lr){1-3} \cmidrule(lr){4-5} \cmidrule(lr){6-12}
500 & d & 10 & 345 & 2 & 239 & 69 & 26 & 5 & 3317 & 220 & 895 \\
500 & d & 25 & 132 & 5 & 74 & 48 & 49 & 3 & 783 & 91 & 108 \\
500 & d & 50 & 9 & 1 & 5 & 80 & 17 & 3 & 88 & 0 & 29 \\
500 & d & 75 & 9 & 1 & 5 & 81 & 16 & 3 & 88 & 0 & 29 \\
500 & d & 100 & 9 & 1 & 5 & 75 & 22 & 3 & 88 & 0 & 29 \\
\cmidrule(lr){1-3} \cmidrule(lr){4-5} \cmidrule(lr){6-12}
500 & w & 10 & 345 & 2 & 260 & 65 & 31 & 4 & 3317 & 220 & 895 \\
500 & w & 25 & 132 & 5 & 81 & 52 & 44 & 4 & 783 & 91 & 108 \\
500 & w & 50 & 9 & 1 & 6 & 80 & 17 & 3 & 88 & 0 & 29 \\
500 & w & 75 & 9 & 1 & 6 & 78 & 19 & 3 & 88 & 0 & 29 \\
500 & w & 100 & 9 & 1 & 8 & 66 & 30 & 4 & 88 & 0 & 29 \\
\hline
\hline
1000 & b & 10 & 714 & 2 & 565 & 73 & 20 & 7 & 13407 & 321 & 5394 \\
1000 & b & 25 & 1548 & 5 & 1030 & 73 & 15 & 12 & 21263 & 62 & 8504 \\
1000 & b & 50 & 726 & 3 & 382 & 57 & 39 & 4 & 8495 & 1348 & 1949 \\
1000 & b & 75 & 84 & 1 & 104 & 35 & 62 & 3 & 1139 & 0 & 471 \\
1000 & b & 100 & 84 & 1 & 46 & 80 & 17 & 3 & 1139 & 0 & 471 \\
\cmidrule(lr){1-3} \cmidrule(lr){4-5} \cmidrule(lr){6-12}
1000 & d & 10 & 942 & 2 & 668 & 70 & 22 & 8 & 17463 & 321 & 7122 \\
1000 & d & 25 & 1716 & 5 & 1000 & 69 & 20 & 11 & 25514 & 62 & 10331 \\
1000 & d & 50 & 1422 & 3 & 628 & 57 & 39 & 4 & 21496 & 2022 & 6832 \\
1000 & d & 75 & 164 & 1 & 75 & 72 & 25 & 3 & 3004 & 0 & 1301 \\
1000 & d & 100 & 164 & 1 & 74 & 73 & 24 & 3 & 3004 & 0 & 1301 \\
\cmidrule(lr){1-3} \cmidrule(lr){4-5} \cmidrule(lr){6-12}
1000 & w & 10 & 966 & 2 & 704 & 68 & 27 & 5 & 26091 & 321 & 11295 \\
1000 & w & 25 & 1782 & 5 & 1278 & 62 & 31 & 7 & 48782 & 62 & 21407 \\
1000 & w & 50 & 1513 & 3 & 923 & 52 & 45 & 3 & 45948 & 2022 & 18589 \\
1000 & w & 75 & 174 & 1 & 104 & 68 & 29 & 3 & 6742 & 0 & 3116 \\
1000 & w & 100 & 174 & 1 & 115 & 62 & 36 & 2 & 6742 & 0 & 3116 \\
\hline
\hline
1500 & b & 10 & 36864 & 2 & 16795 & 39 & 30 & 31 & 154653 & 1489 & 30324 \\
1500 & b & 25 & 40872 & 8 & 12800 & 48 & 41 & 11 & 161295 & 942 & 32380 \\
1500 & b & 50 & 11240 & 3 & 3257 & 40 & 47 & 13 & 44319 & 214 & 9010 \\
1500 & b & 75 & 9832 & 2 & 2560 & 26 & 54 & 20 & 34553 & 2082 & 4534 \\
1500 & b & 100 & 1004 & 1 & 335 & 46 & 39 & 15 & 4109 & 0 & 902 \\
\cmidrule(lr){1-3} \cmidrule(lr){4-5} \cmidrule(lr){6-12}
1500 & d & 10 & 44630 & 4 & 24290 & 40 & 32 & 28 & 187349 & 1109 & 37230 \\
1500 & d & 25 & 44032 & 8 & 14717 & 48 & 39 & 13 & 163057 & 942 & 29409 \\
1500 & d & 50 & 20216 & 4 & 5734 & 43 & 44 & 13 & 77815 & 321 & 14975 \\
1500 & d & 75 & 19516 & 2 & 5903 & 24 & 56 & 20 & 66685 & 4450 & 7864 \\
1500 & d & 100 & 1832 & 2 & 854 & 43 & 40 & 17 & 7199 & 0 & 1468 \\
\cmidrule(lr){1-3} \cmidrule(lr){4-5} \cmidrule(lr){6-12}
1500 & w & 10 & 43938 & 4 & 20593 & 40 & 32 & 28 & 182194 & 805 & 35633 \\
1500 & w & 25 & 45728 & 8 & 16103 & 48 & 42 & 10 & 173389 & 942 & 30894 \\
1500 & w & 50 & 20920 & 4 & 6712 & 46 & 42 & 12 & 82577 & 321 & 15757 \\
1500 & w & 75 & 22016 & 2 & 5567 & 29 & 49 & 22 & 74648 & 5083 & 8555 \\
1500 & w & 100 & 1924 & 2 & 663 & 46 & 39 & 15 & 7610 & 0 & 1501 \\

\hline

\end{tabular}}
\caption{Summary of the paths found for each instance}
\label{tab:summary}
\end{table}

\begingroup
\setlength{\tabcolsep}{4pt}  
\renewcommand{\arraystretch}{0.95} 

\begin{table}[ht]
    \centering
\resizebox{0.9\columnwidth}{!}{%
\begin{tabular}{ccrcrrrrrrccccccc}	
	
    \multicolumn{3}{c}{Instance} & \multicolumn{2}{c}{}  & \multicolumn{2}{c}{Fuel} & \multicolumn{3}{c}{Time (h)} & \multicolumn{7}{c}{Number of stops}\\
    \cmidrule(lr){1-3} \cmidrule(lr){6-7} \cmidrule(lr){8-10} \cmidrule(lr){11-17}

    km & HoS & fuel & \begin{tabular}[r]{@{}c@{}}Path \\ \# \end{tabular} & \begin{tabular}[r]{@{}c@{}}Distance \\ (km) \end{tabular} & \begin{tabular}[r]{@{}c@{}}Cost \\ (\euro) \end{tabular} & Liters & Total & Travel & Stop & {\footnotesize F} & {\footnotesize FB} & {\footnotesize FD} & {\footnotesize FW} & {\footnotesize B} & {\footnotesize D} & {\footnotesize W} \\
    \hline \hline
    
500 & b & 10 & CP & 516.5 & 207.0 & 147.6 & 16.42 & 5.17 & 11.25 & 1 & 0 & 0 & 0 & 0 & 1 & 0 \\
&  &  & 1 & 468.8 & 148.4 & 133.9 & 15.69 & 4.69 & 11.00 & 0 & 0 & 1 & 0 & 0 & 0 & 0 \\
&  &  & 2 & 470.9 & 145.0 & 134.5 & 15.96 & 4.71 & 11.25 & 1 & 0 & 0 & 0 & 0 & 1 & 0 \\
\cmidrule(lr){4-17}
500 & b & 25 & CP & 518.5 & 200.9 & 148.1 & 16.18 & 5.18 & 11.00 & 0 & 0 & 1 & 0 & 0 & 0 & 0 \\
&  &  & 1 & 471.7 & 182.7 & 134.8 & 15.72 & 4.72 & 11.00 & 0 & 0 & 1 & 0 & 0 & 0 & 0 \\
&  &  & 2 & 472.8 & 175.6 & 135.1 & 15.73 & 4.73 & 11.00 & 0 & 0 & 1 & 0 & 0 & 0 & 0 \\
&  &  & 3 & 475.7 & 172.6 & 135.9 & 15.76 & 4.76 & 11.00 & 0 & 0 & 1 & 0 & 0 & 0 & 0 \\
&  &  & 4 & 488.2 & 172.1 & 139.5 & 15.88 & 4.88 & 11.00 & 0 & 0 & 1 & 0 & 0 & 0 & 0 \\
&  &  & 5 & 489.1 & 171.1 & 139.7 & 15.89 & 4.89 & 11.00 & 0 & 0 & 1 & 0 & 0 & 0 & 0 \\
\cmidrule(lr){4-17}
500 & b & 50 & CP & 516.5 & 206.6 & 147.6 & 16.17 & 5.17 & 11.00 & 0 & 0 & 0 & 0 & 0 & 1 & 0 \\
&  &  & 1 & 472.3 & 188.9 & 134.9 & 15.72 & 4.72 & 11.00 & 0 & 0 & 0 & 0 & 0 & 1 & 0 \\
\cmidrule(lr){4-17}
500 & b & 75 & CP & 516.5 & 206.6 & 147.6 & 16.17 & 5.17 & 11.00 & 0 & 0 & 0 & 0 & 0 & 1 & 0 \\
&  &  & 1 & 472.3 & 188.9 & 134.9 & 15.72 & 4.72 & 11.00 & 0 & 0 & 0 & 0 & 0 & 1 & 0 \\
\cmidrule(lr){4-17}
500 & b & 100 & CP & 516.5 & 206.6 & 147.6 & 16.17 & 5.17 & 11.00 & 0 & 0 & 0 & 0 & 0 & 1 & 0 \\
&  &  & 1 & 472.3 & 188.9 & 134.9 & 15.72 & 4.72 & 11.00 & 0 & 0 & 0 & 0 & 0 & 1 & 0 \\
\hline \hline
500 & d & 10 & CP & 516.5 & 207.0 & 147.6 & 6.17 & 5.17 & 1.00 & 1 & 0 & 0 & 0 & 1 & 0 & 0 \\
&  &  & 1 & 468.8 & 148.4 & 133.9 & 5.71 & 4.71 & 0.75 & 0 & 1 & 0 & 0 & 0 & 0 & 0 \\
&  &  & 2 & 470.9 & 145.0 & 134.5 & 5.44 & 4.69 & 1.00 & 1 & 0 & 0 & 0 & 1 & 0 & 0 \\
\cmidrule(lr){4-17}
500 & d & 25 & CP & 518.5 & 200.9 & 148.1 & 5.93 & 5.18 & 0.75 & 0 & 1 & 0 & 0 & 0 & 0 & 0 \\
&  &  & 1 & 471.7 & 182.7 & 134.8 & 5.47 & 4.72 & 0.75 & 0 & 1 & 0 & 0 & 0 & 0 & 0 \\
&  &  & 2 & 472.8 & 175.6 & 135.1 & 5.48 & 4.73 & 0.75 & 0 & 1 & 0 & 0 & 0 & 0 & 0 \\
&  &  & 3 & 475.7 & 172.6 & 135.9 & 5.51 & 4.76 & 0.75 & 0 & 1 & 0 & 0 & 0 & 0 & 0 \\
&  &  & 4 & 488.2 & 172.1 & 139.5 & 5.63 & 4.88 & 0.75 & 0 & 1 & 0 & 0 & 0 & 0 & 0 \\
&  &  & 5 & 489.1 & 171.1 & 139.7 & 5.64 & 4.89 & 0.75 & 0 & 1 & 0 & 0 & 0 & 0 & 0 \\
\cmidrule(lr){4-17}
500 & d & 50 & CP & 516.5 & 206.6 & 147.6 & 5.92 & 5.17 & 0.75 & 0 & 0 & 0 & 0 & 1 & 0 & 0 \\
&  &  & 1 & 472.3 & 188.9 & 134.9 & 5.47 & 4.72 & 0.75 & 0 & 0 & 0 & 0 & 1 & 0 & 0 \\
\cmidrule(lr){4-17}
500 & d & 75 & CP & 516.5 & 206.6 & 147.6 & 5.92 & 5.17 & 0.75 & 0 & 0 & 0 & 0 & 1 & 0 & 0 \\
&  &  & 1 & 472.3 & 188.9 & 134.9 & 5.47 & 4.72 & 0.75 & 0 & 0 & 0 & 0 & 1 & 0 & 0 \\
\cmidrule(lr){4-17}
500 & d & 100 & CP & 516.5 & 206.6 & 147.6 & 5.92 & 5.17 & 0.75 & 0 & 0 & 0 & 0 & 1 & 0 & 0 \\
&  &  & 1 & 472.3 & 188.9 & 134.9 & 5.47 & 4.72 & 0.75 & 0 & 0 & 0 & 0 & 1 & 0 & 0 \\
\hline \hline
500 & w & 10 & CP & 516.5 & 207.0 & 147.6 & 6.17 & 5.17 & 1.00 & 1 & 0 & 0 & 0 & 1 & 0 & 0 \\
&  &  & 1 & 468.8 & 148.4 & 133.9 & 5.44 & 4.69 & 0.75 & 0 & 1 & 0 & 0 & 0 & 0 & 0 \\
&  &  & 2 & 470.9 & 145.0 & 134.5 & 5.71 & 4.71 & 1.00 & 1 & 0 & 0 & 0 & 1 & 0 & 0 \\
\cmidrule(lr){4-17}
500 & w & 25 & CP & 518.5 & 200.9 & 148.1 & 5.93 & 5.18 & 0.75 & 0 & 1 & 0 & 0 & 0 & 0 & 0 \\
&  &  & 1 & 471.7 & 182.7 & 134.8 & 5.47 & 4.72 & 0.75 & 0 & 1 & 0 & 0 & 0 & 0 & 0 \\
&  &  & 2 & 472.8 & 175.6 & 135.1 & 5.48 & 4.73 & 0.75 & 0 & 1 & 0 & 0 & 0 & 0 & 0 \\
&  &  & 3 & 475.7 & 172.6 & 135.9 & 5.51 & 4.76 & 0.75 & 0 & 1 & 0 & 0 & 0 & 0 & 0 \\
&  &  & 4 & 488.2 & 172.1 & 139.5 & 5.63 & 4.88 & 0.75 & 0 & 1 & 0 & 0 & 0 & 0 & 0 \\
&  &  & 5 & 489.1 & 171.1 & 139.7 & 5.64 & 4.89 & 0.75 & 0 & 1 & 0 & 0 & 0 & 0 & 0 \\
\cmidrule(lr){4-17}
500 & w & 50 & CP & 516.5 & 206.6 & 147.6 & 5.92 & 5.17 & 0.75 & 0 & 0 & 0 & 0 & 1 & 0 & 0 \\
&  &  & 1 & 472.3 & 188.9 & 134.9 & 5.47 & 4.72 & 0.75 & 0 & 0 & 0 & 0 & 1 & 0 & 0 \\
\cmidrule(lr){4-17}
500 & w & 75 & CP & 516.5 & 206.6 & 147.6 & 5.92 & 5.17 & 0.75 & 0 & 0 & 0 & 0 & 1 & 0 & 0 \\
&  &  & 1 & 472.3 & 188.9 & 134.9 & 5.47 & 4.72 & 0.75 & 0 & 0 & 0 & 0 & 1 & 0 & 0 \\
\cmidrule(lr){4-17}
500 & w & 100 & CP & 516.5 & 206.6 & 147.6 & 5.92 & 5.17 & 0.75 & 0 & 0 & 0 & 0 & 1 & 0 & 0 \\
&  &  & 1 & 472.3 & 188.9 & 134.9 & 5.47 & 4.72 & 0.75 & 0 & 0 & 0 & 0 & 1 & 0 & 0 \\

    \hline 
\end{tabular}}
\caption{Details of the non-dominated paths found for the instances of 500 km}
\label{tab:results}
\end{table}

\begin{table}[ht]
	\centering
	\resizebox{0.9\columnwidth}{!}{%
		\begin{tabular}{ccrcrrrrrrccccccc}
			
    	\multicolumn{3}{c}{Instance} & \multicolumn{2}{c}{}  & \multicolumn{2}{c}{Fuel} & \multicolumn{3}{c}{Time (h)} & \multicolumn{7}{c}{Number of stops}\\
		\cmidrule(lr){1-3} \cmidrule(lr){6-7} \cmidrule(lr){8-10} \cmidrule(lr){11-17}
			
    km & HoS & fuel & \begin{tabular}[r]{@{}c@{}}Path \\ \# \end{tabular} & \begin{tabular}[r]{@{}c@{}}Distance \\ (km) \end{tabular} & \begin{tabular}[r]{@{}c@{}}Cost \\ (\euro) \end{tabular} & Liters & Total & Travel & Stop & {\footnotesize F} & {\footnotesize FB} & {\footnotesize FD} & {\footnotesize FW} & {\footnotesize B} & {\footnotesize D} & {\footnotesize W} \\
			\hline \hline
		
1000 & b & 10 & CP & 930.2 & 388.7 & 265.8 & 21.30 & 9.30 & 12.00 & 1 & 0 & 0 & 0 & 1 & 1 & 0 \\
&  &  & 1 & 899.8 & 322.5 & 257.1 & 21.50 & 9.00 & 12.50 & 0 & 0 & 1 & 0 & 2 & 0 & 0 \\
&  &  & 2 & 911.9 & 326.9 & 260.5 & 21.12 & 9.12 & 12.00 & 1 & 0 & 0 & 0 & 1 & 1 & 0 \\
\cmidrule(lr){4-17}
1000 & b & 25 & CP & 929.8 & 411.7 & 265.7 & 21.05 & 9.30 & 11.75 & 0 & 0 & 1 & 0 & 1 & 0 & 0 \\
&  &  & 1 & 929.7 & 392.8 & 265.6 & 21.05 & 9.30 & 11.75 & 0 & 0 & 1 & 0 & 1 & 0 & 0 \\
&  &  & 2 & 929.8 & 346.1 & 265.6 & 21.30 & 9.30 & 12.00 & 1 & 0 & 0 & 0 & 1 & 1 & 0 \\
&  &  & 3 & 931.7 & 372.7 & 266.2 & 21.07 & 9.32 & 11.75 & 0 & 0 & 1 & 0 & 1 & 0 & 0 \\
&  &  & 4 & 933.4 & 352.6 & 266.7 & 21.08 & 9.33 & 11.75 & 0 & 0 & 1 & 0 & 1 & 0 & 0 \\
&  &  & 5 & 936.5 & 348.5 & 267.6 & 21.11 & 9.36 & 11.75 & 0 & 0 & 1 & 0 & 1 & 0 & 0 \\
\cmidrule(lr){4-17}
1000 & b & 50 & CP & 927.0 & 373.6 & 264.8 & 21.02 & 9.27 & 11.75 & 0 & 1 & 0 & 0 & 0 & 1 & 0 \\
&  &  & 1 & 926.1 & 371.6 & 264.6 & 21.01 & 9.26 & 11.75 & 0 & 1 & 0 & 0 & 0 & 1 & 0 \\
&  &  & 2 & 927.1 & 365.4 & 264.9 & 21.02 & 9.27 & 11.75 & 0 & 1 & 0 & 0 & 0 & 1 & 0 \\
&  &  & 3 & 928.8 & 363.4 & 265.4 & 21.04 & 9.29 & 11.75 & 0 & 1 & 0 & 0 & 0 & 1 & 0 \\
\cmidrule(lr){4-17}
1000 & b & 75 & CP & 927.0 & 370.8 & 264.8 & 21.02 & 9.27 & 11.75 & 0 & 0 & 0 & 0 & 1 & 1 & 0 \\
&  &  & 1 & 928.9 & 371.6 & 265.4 & 21.04 & 9.29 & 11.75 & 0 & 0 & 0 & 0 & 1 & 1 & 0 \\
\cmidrule(lr){4-17}
1000 & b & 100 & CP & 927.0 & 370.8 & 264.8 & 21.02 & 9.27 & 11.75 & 0 & 0 & 0 & 0 & 1 & 1 & 0 \\
&  &  & 1 & 928.9 & 371.6 & 265.4 & 21.04 & 9.29 & 11.75 & 0 & 0 & 0 & 0 & 1 & 1 & 0 \\
\hline \hline
1000 & d & 10 & CP & 930.2 & 388.7 & 265.8 & 22.05 & 9.30 & 12.75 & 1 & 0 & 0 & 0 & 2 & 1 & 0 \\
&  &  & 1 & 899.8 & 322.5 & 257.1 & 21.50 & 9.00 & 12.50 & 0 & 0 & 1 & 0 & 2 & 0 & 0 \\
&  &  & 2 & 911.9 & 326.9 & 260.5 & 21.12 & 9.12 & 12.00 & 1 & 0 & 0 & 0 & 1 & 1 & 0 \\
\cmidrule(lr){4-17}
1000 & d & 25 & CP & 929.8 & 411.7 & 265.7 & 21.80 & 9.30 & 12.50 & 0 & 1 & 0 & 0 & 1 & 1 & 0 \\
&  &  & 1 & 929.7 & 392.8 & 265.6 & 21.05 & 9.30 & 11.75 & 0 & 0 & 1 & 0 & 1 & 0 & 0 \\
&  &  & 2 & 929.8 & 346.1 & 265.6 & 21.30 & 9.30 & 12.00 & 1 & 0 & 0 & 0 & 1 & 1 & 0 \\
&  &  & 3 & 931.7 & 372.7 & 266.2 & 21.07 & 9.32 & 11.75 & 0 & 0 & 1 & 0 & 1 & 0 & 0 \\
&  &  & 4 & 933.4 & 352.6 & 266.7 & 21.08 & 9.33 & 11.75 & 0 & 0 & 1 & 0 & 1 & 0 & 0 \\
&  &  & 5 & 936.5 & 348.5 & 267.6 & 21.11 & 9.36 & 11.75 & 0 & 0 & 1 & 0 & 1 & 0 & 0 \\
\cmidrule(lr){4-17}
1000 & d & 50 & CP & 927.0 & 373.6 & 264.8 & 21.02 & 9.27 & 11.75 & 0 & 0 & 1 & 0 & 1 & 0 & 0 \\
&  &  & 1 & 926.1 & 371.6 & 264.6 & 21.01 & 9.26 & 11.75 & 0 & 1 & 0 & 0 & 0 & 1 & 0 \\
&  &  & 2 & 927.1 & 365.4 & 264.9 & 21.02 & 9.27 & 11.75 & 0 & 1 & 0 & 0 & 0 & 1 & 0 \\
&  &  & 3 & 928.8 & 363.4 & 265.4 & 21.04 & 9.29 & 11.75 & 0 & 1 & 0 & 0 & 0 & 1 & 0 \\
\cmidrule(lr){4-17}
1000 & d & 75 & CP & 927.0 & 370.8 & 264.8 & 21.77 & 9.27 & 12.50 & 0 & 0 & 0 & 0 & 2 & 1 & 0 \\
&  &  & 1 & 928.9 & 371.6 & 265.4 & 21.04 & 9.29 & 11.75 & 0 & 0 & 0 & 0 & 1 & 1 & 0 \\
\cmidrule(lr){4-17}
1000 & d & 100 & CP & 927.0 & 370.8 & 264.8 & 21.77 & 9.27 & 12.50 & 0 & 0 & 0 & 0 & 2 & 1 & 0 \\
&  &  & 1 & 928.9 & 371.6 & 265.4 & 21.04 & 9.29 & 11.75 & 0 & 0 & 0 & 0 & 1 & 1 & 0 \\
\hline \hline
1000 & w & 10 & CP & 930.2 & 388.7 & 265.8 & 56.05 & 9.30 & 46.75 & 1 & 0 & 0 & 0 & 2 & 0 & 1 \\
&  &  & 1 & 899.8 & 322.5 & 257.1 & 55.50 & 9.00 & 46.50 & 0 & 0 & 0 & 1 & 2 & 0 & 0 \\
&  &  & 2 & 911.9 & 326.9 & 260.5 & 55.12 & 9.12 & 46.00 & 1 & 0 & 0 & 0 & 1 & 0 & 1 \\
\cmidrule(lr){4-17}
1000 & w & 25 & CP & 929.8 & 411.7 & 265.7 & 55.80 & 9.30 & 46.50 & 0 & 1 & 0 & 0 & 1 & 0 & 1 \\
&  &  & 1 & 929.7 & 392.8 & 265.6 & 55.05 & 9.30 & 45.75 & 0 & 0 & 0 & 1 & 1 & 0 & 0 \\
&  &  & 2 & 929.8 & 346.1 & 265.6 & 55.30 & 9.30 & 46.00 & 1 & 0 & 0 & 0 & 1 & 0 & 1 \\
&  &  & 3 & 931.7 & 372.7 & 266.2 & 55.07 & 9.32 & 45.75 & 0 & 0 & 0 & 1 & 1 & 0 & 0 \\
&  &  & 4 & 933.4 & 352.6 & 266.7 & 55.08 & 9.33 & 45.75 & 0 & 0 & 0 & 1 & 1 & 0 & 0 \\
&  &  & 5 & 936.5 & 348.5 & 267.6 & 55.11 & 9.36 & 45.75 & 0 & 0 & 0 & 1 & 1 & 0 & 0 \\
\cmidrule(lr){4-17}
1000 & w & 50 & CP & 927.0 & 373.6 & 264.8 & 55.02 & 9.27 & 45.75 & 0 & 0 & 0 & 1 & 1 & 0 & 0 \\
&  &  & 1 & 926.1 & 371.6 & 264.6 & 55.01 & 9.26 & 45.75 & 0 & 1 & 0 & 0 & 0 & 0 & 1 \\
&  &  & 2 & 927.1 & 365.4 & 264.9 & 55.02 & 9.27 & 45.75 & 0 & 1 & 0 & 0 & 0 & 0 & 1 \\
&  &  & 3 & 928.8 & 363.4 & 265.4 & 55.04 & 9.29 & 45.75 & 0 & 1 & 0 & 0 & 0 & 0 & 1 \\
\cmidrule(lr){4-17}
1000 & w & 75 & CP & 927.0 & 370.8 & 264.8 & 55.77 & 9.27 & 46.50 & 0 & 0 & 0 & 0 & 2 & 0 & 1 \\
&  &  & 1 & 928.9 & 371.6 & 265.4 & 55.04 & 9.29 & 45.75 & 0 & 0 & 0 & 0 & 1 & 0 & 1 \\
\cmidrule(lr){4-17}
1000 & w & 100 & CP & 927.0 & 370.8 & 264.8 & 55.77 & 9.27 & 46.50 & 0 & 0 & 0 & 0 & 2 & 0 & 1 \\
&  &  & 1 & 928.9 & 371.6 & 265.4 & 55.04 & 9.29 & 45.75 & 0 & 0 & 0 & 0 & 1 & 0 & 1 \\
			
			\hline 
	\end{tabular}}
	\caption{Details of the non-dominated paths found for the instances of 1000 km}
	\label{tab:results1000}
\end{table}

\begin{table}[ht]
	\centering
	\resizebox{0.9\columnwidth}{!}{%
		\begin{tabular}{ccrcrrrrrrccccccc}
			\multicolumn{3}{c}{Instance} & \multicolumn{2}{c}{}  & \multicolumn{2}{c}{Fuel} & \multicolumn{3}{c}{Time (h)} & \multicolumn{7}{c}{Number of stops}\\
			\cmidrule(lr){1-3} \cmidrule(lr){6-7} \cmidrule(lr){8-10} \cmidrule(lr){11-17}
			
    km & HoS & fuel & \begin{tabular}[r]{@{}c@{}}Path \\ \# \end{tabular} & \begin{tabular}[r]{@{}c@{}}Distance \\ (km) \end{tabular} & \begin{tabular}[r]{@{}c@{}}Cost \\ (\euro) \end{tabular} & Liters & Total & Travel & Stop & {\footnotesize F} & {\footnotesize FB} & {\footnotesize FD} & {\footnotesize FW} & {\footnotesize B} & {\footnotesize D} & {\footnotesize W} \\
			\hline \hline

1500 & b & 10 & CP & 1508.9 & 640.2 & 431.1 & 38.84 & 15.09 & 23.75 & 1 & 0 & 0 & 0 & 2 & 2 & 0 \\
&  &  & 1 & 1445.3 & 525.4 & 413.0 & 37.95 & 14.45 & 23.50 & 0 & 0 & 1 & 0 & 2 & 1 & 0 \\
&  &  & 2 & 1448.4 & 526.5 & 413.8 & 37.48 & 14.48 & 23.00 & 1 & 0 & 0 & 0 & 1 & 2 & 0 \\
\cmidrule(lr){4-17}
1500 & b & 25 & CP & 1505.8 & 560.1 & 430.2 & 38.56 & 15.06 & 23.50 & 0 & 0 & 1 & 0 & 2 & 1 & 0 \\
&  &  & 1 & 1447.2 & 627.5 & 413.5 & 37.22 & 14.47 & 22.75 & 0 & 0 & 1 & 0 & 1 & 1 & 0 \\
&  &  & 2 & 1448.0 & 589.9 & 413.7 & 37.23 & 14.48 & 22.75 & 0 & 0 & 1 & 0 & 1 & 1 & 0 \\
&  &  & 3 & 1451.0 & 539.7 & 414.6 & 37.51 & 14.51 & 23.00 & 1 & 0 & 0 & 0 & 1 & 2 & 0 \\
&  &  & 4 & 1452.8 & 537.1 & 415.1 & 37.53 & 14.53 & 23.00 & 1 & 0 & 0 & 0 & 1 & 2 & 0 \\
&  &  & 5 & 1464.3 & 544.5 & 418.4 & 37.39 & 14.64 & 22.75 & 0 & 0 & 1 & 0 & 1 & 1 & 0 \\
&  &  & 6 & 1466.0 & 541.9 & 418.9 & 37.41 & 14.66 & 22.75 & 0 & 0 & 1 & 0 & 1 & 1 & 0 \\
&  &  & 7 & 1478.9 & 489.5 & 422.5 & 37.54 & 14.79 & 22.75 & 0 & 0 & 1 & 0 & 1 & 1 & 0 \\
&  &  & 8 & 1490.7 & 485.0 & 425.9 & 37.66 & 14.91 & 22.75 & 0 & 0 & 1 & 0 & 1 & 1 & 0 \\
\cmidrule(lr){4-17}
1500 & b & 50 & CP & 1508.9 & 574.9 & 431.1 & 38.59 & 15.09 & 23.50 & 0 & 1 & 0 & 0 & 1 & 2 & 0 \\
&  &  & 1 & 1412.6 & 547.9 & 403.6 & 37.13 & 14.13 & 23.00 & 1 & 0 & 0 & 0 & 1 & 2 & 0 \\
&  &  & 2 & 1427.8 & 512.9 & 407.9 & 37.78 & 14.28 & 23.50 & 0 & 1 & 0 & 0 & 1 & 2 & 0 \\
&  &  & 3 & 1445.9 & 518.5 & 413.1 & 37.21 & 14.46 & 22.75 & 0 & 0 & 1 & 0 & 1 & 1 & 0 \\
\cmidrule(lr){4-17}
1500 & b & 75 & CP & 1507.5 & 593.6 & 430.7 & 37.82 & 15.07 & 22.75 & 0 & 0 & 1 & 0 & 1 & 1 & 0 \\
&  &  & 1 & 1448.4 & 565.2 & 413.8 & 37.23 & 14.48 & 22.75 & 0 & 1 & 0 & 0 & 0 & 2 & 0 \\
&  &  & 2 & 1477.1 & 560.2 & 422.0 & 37.52 & 14.77 & 22.75 & 0 & 1 & 0 & 0 & 0 & 2 & 0 \\
\cmidrule(lr){4-17}
1500 & b & 100 & CP & 1508.9 & 603.6 & 431.1 & 38.59 & 15.09 & 23.50 & 0 & 0 & 0 & 0 & 2 & 2 & 0 \\
&  &  & 1 & 1448.4 & 579.4 & 413.8 & 37.23 & 14.48 & 22.75 & 0 & 0 & 0 & 0 & 1 & 2 & 0 \\	


			
			\hline \hline

1500 & d & 10 & CP & 1495.1 & 634.3 & 427.2 & 28.45 & 14.95 & 13.50 & 1 & 0 & 0 & 0 & 3 & 1 & 0 \\
&  &  & 1 & 1445.0 & 525.3 & 412.9 & 38.70 & 14.45 & 24.25 & 0 & 0 & 1 & 0 & 3 & 1 & 0 \\
&  &  & 2 & 1445.3 & 525.4 & 413.0 & 37.95 & 14.45 & 23.50 & 0 & 0 & 1 & 0 & 2 & 1 & 0 \\
&  &  & 3 & 1446.7 & 526.0 & 413.4 & 27.72 & 14.47 & 13.25 & 0 & 1 & 0 & 0 & 2 & 1 & 0 \\
&  &  & 4 & 1448.4 & 526.5 & 413.8 & 27.23 & 14.48 & 12.75 & 1 & 0 & 0 & 0 & 2 & 1 & 0 \\
\cmidrule(lr){4-17}
1500 & d & 25 & CP & 1495.4 & 556.3 & 427.3 & 28.20 & 14.95 & 13.25 & 0 & 1 & 0 & 0 & 2 & 1 & 0 \\
&  &  & 1 & 1447.2 & 627.5 & 413.5 & 26.97 & 14.47 & 12.50 & 0 & 1 & 0 & 0 & 1 & 1 & 0 \\
&  &  & 2 & 1448.0 & 589.9 & 413.7 & 26.98 & 14.48 & 12.50 & 0 & 1 & 0 & 0 & 1 & 1 & 0 \\
&  &  & 3 & 1451.0 & 539.7 & 414.6 & 27.26 & 14.51 & 12.75 & 1 & 0 & 0 & 0 & 2 & 1 & 0 \\
&  &  & 4 & 1452.8 & 537.1 & 415.1 & 27.28 & 14.53 & 12.75 & 1 & 0 & 0 & 0 & 2 & 1 & 0 \\
&  &  & 5 & 1464.3 & 544.5 & 418.4 & 27.14 & 14.64 & 12.50 & 0 & 1 & 0 & 0 & 1 & 1 & 0 \\
&  &  & 6 & 1466.0 & 541.9 & 418.9 & 27.16 & 14.66 & 12.50 & 0 & 1 & 0 & 0 & 1 & 1 & 0 \\
&  &  & 7 & 1478.9 & 489.5 & 422.5 & 27.29 & 14.79 & 12.50 & 0 & 1 & 0 & 0 & 1 & 1 & 0 \\
&  &  & 8 & 1490.7 & 485.0 & 425.9 & 27.41 & 14.91 & 12.50 & 0 & 1 & 0 & 0 & 1 & 1 & 0 \\
\cmidrule(lr){4-17}
1500 & d & 50 & CP & 1507.5 & 574.4 & 430.7 & 27.57 & 15.07 & 12.50 & 0 & 0 & 1 & 0 & 2 & 0 & 0 \\
&  &  & 1 & 1412.6 & 547.9 & 403.6 & 26.88 & 14.13 & 12.75 & 1 & 0 & 0 & 0 & 2 & 1 & 0 \\
&  &  & 2 & 1427.8 & 512.9 & 407.9 & 37.78 & 14.28 & 23.50 & 0 & 1 & 0 & 0 & 1 & 2 & 0 \\
&  &  & 3 & 1445.9 & 518.5 & 413.1 & 26.96 & 14.46 & 12.50 & 0 & 0 & 1 & 0 & 2 & 0 & 0 \\
&  &  & 4 & 1498.3 & 517.8 & 428.1 & 27.48 & 14.98 & 12.50 & 0 & 0 & 1 & 0 & 2 & 0 & 0 \\
\cmidrule(lr){4-17}
1500 & d & 75 & CP & 1495.1 & 598.7 & 427.2 & 28.20 & 14.95 & 13.25 & 0 & 1 & 0 & 0 & 2 & 1 & 0 \\
&  &  & 1 & 1448.4 & 565.2 & 413.8 & 26.98 & 14.48 & 12.50 & 0 & 1 & 0 & 0 & 1 & 1 & 0 \\
&  &  & 2 & 1477.1 & 560.2 & 422.0 & 27.27 & 14.77 & 12.50 & 0 & 1 & 0 & 0 & 1 & 1 & 0 \\
\cmidrule(lr){4-17}
1500 & d & 100 & CP & 1495.1 & 598.0 & 427.2 & 28.20 & 14.95 & 13.25 & 0 & 0 & 0 & 0 & 3 & 1 & 0 \\
&  &  & 1 & 1447.2 & 578.9 & 413.5 & 37.97 & 14.47 & 23.50 & 0 & 0 & 0 & 0 & 2 & 2 & 0 \\
&  &  & 2 & 1448.4 & 579.4 & 413.8 & 26.98 & 14.48 & 12.50 & 0 & 0 & 0 & 0 & 2 & 1 & 0 \\

			\hline			
	\end{tabular}}
	\caption{Details of the non-dominated paths found for the instances of 1500 km with b and d}
	\label{tab:results1500_2}
\end{table}

\begin{table}[ht]
	\centering
	\resizebox{0.9\columnwidth}{!}{%
		\begin{tabular}{ccrcrrrrrrccccccc}
			\multicolumn{3}{c}{Instance} & \multicolumn{2}{c}{}  & \multicolumn{2}{c}{Fuel} & \multicolumn{3}{c}{Time (h)} & \multicolumn{7}{c}{Number of stops}\\
			\cmidrule(lr){1-3} \cmidrule(lr){6-7} \cmidrule(lr){8-10} \cmidrule(lr){11-17}
			
			km & HoS & fuel & \begin{tabular}[r]{@{}c@{}}Path \\ \# \end{tabular} & \begin{tabular}[r]{@{}c@{}}Distance \\ (km) \end{tabular} & \begin{tabular}[r]{@{}c@{}}Cost \\ (\euro) \end{tabular} & Liters & Total & Travel & Stop & {\footnotesize F} & {\footnotesize FB} & {\footnotesize FD} & {\footnotesize FW} & {\footnotesize B} & {\footnotesize D} & {\footnotesize W} \\
			\hline \hline

1500 & w & 10 & CP & 1495.1 & 634.3 & 427.2 & 62.45 & 14.95 & 47.50 & 1 & 0 & 0 & 0 & 3 & 0 & 1 \\
&  &  & 1 & 1445.0 & 525.3 & 412.9 & 72.70 & 14.45 & 58.25 & 0 & 0 & 0 & 1 & 3 & 1 & 0 \\
&  &  & 2 & 1445.3 & 525.4 & 413.0 & 71.95 & 14.45 & 57.50 & 0 & 0 & 0 & 1 & 2 & 1 & 0 \\
&  &  & 3 & 1446.7 & 526.0 & 413.4 & 61.72 & 14.47 & 47.25 & 0 & 1 & 0 & 0 & 2 & 0 & 1 \\
&  &  & 4 & 1448.4 & 526.5 & 413.8 & 61.23 & 14.48 & 46.75 & 1 & 0 & 0 & 0 & 2 & 0 & 1 \\
\cmidrule(lr){4-17}
1500 & w & 25 & CP & 1495.4 & 556.3 & 427.3 & 62.20 & 14.95 & 47.25 & 0 & 1 & 0 & 0 & 2 & 0 & 1 \\
&  &  & 1 & 1447.2 & 627.5 & 413.5 & 60.97 & 14.47 & 46.50 & 0 & 1 & 0 & 0 & 1 & 0 & 1 \\
&  &  & 2 & 1448.0 & 589.9 & 413.7 & 60.98 & 14.48 & 46.50 & 0 & 1 & 0 & 0 & 1 & 0 & 1 \\
&  &  & 3 & 1451.0 & 539.7 & 414.6 & 61.26 & 14.51 & 46.75 & 1 & 0 & 0 & 0 & 2 & 0 & 1 \\
&  &  & 4 & 1452.8 & 537.1 & 415.1 & 61.28 & 14.53 & 46.75 & 1 & 0 & 0 & 0 & 2 & 0 & 1 \\
&  &  & 5 & 1464.3 & 544.5 & 418.4 & 61.14 & 14.64 & 46.50 & 0 & 1 & 0 & 0 & 1 & 0 & 1 \\
&  &  & 6 & 1466.0 & 541.9 & 418.9 & 61.16 & 14.66 & 46.50 & 0 & 1 & 0 & 0 & 1 & 0 & 1 \\
&  &  & 7 & 1478.9 & 489.5 & 422.5 & 61.29 & 14.79 & 46.50 & 0 & 1 & 0 & 0 & 1 & 0 & 1 \\
&  &  & 8 & 1490.7 & 485.0 & 425.9 & 61.41 & 14.91 & 46.50 & 0 & 1 & 0 & 0 & 1 & 0 & 1 \\
\cmidrule(lr){4-17}
1500 & w & 50 & CP & 1507.5 & 574.4 & 430.7 & 61.57 & 15.07 & 46.50 & 0 & 0 & 0 & 1 & 2 & 0 & 0 \\
&  &  & 1 & 1412.6 & 547.9 & 403.6 & 60.88 & 14.13 & 46.75 & 1 & 0 & 0 & 0 & 2 & 0 & 1 \\
&  &  & 2 & 1427.8 & 512.9 & 407.9 & 71.78 & 14.28 & 57.50 & 0 & 1 & 0 & 0 & 1 & 1 & 1 \\
&  &  & 3 & 1445.9 & 518.5 & 413.1 & 60.96 & 14.46 & 46.50 & 0 & 0 & 0 & 1 & 2 & 0 & 0 \\
&  &  & 4 & 1498.3 & 517.8 & 428.1 & 61.48 & 14.98 & 46.50 & 0 & 0 & 0 & 1 & 2 & 0 & 0 \\
\cmidrule(lr){4-17}
1500 & w & 75 & CP & 1495.1 & 598.7 & 427.2 & 62.20 & 14.95 & 47.25 & 0 & 1 & 0 & 0 & 2 & 0 & 1 \\
&  &  & 1 & 1448.4 & 565.2 & 413.8 & 60.98 & 14.48 & 46.50 & 0 & 1 & 0 & 0 & 1 & 0 & 1 \\
&  &  & 2 & 1477.1 & 560.2 & 422.0 & 61.27 & 14.77 & 46.50 & 0 & 1 & 0 & 0 & 1 & 0 & 1 \\
\cmidrule(lr){4-17}
1500 & w & 100 & CP & 1495.1 & 598.0 & 427.2 & 62.20 & 14.95 & 47.25 & 0 & 0 & 0 & 0 & 3 & 0 & 1 \\
&  &  & 1 & 1447.2 & 578.9 & 413.5 & 71.97 & 14.47 & 57.50 & 0 & 0 & 0 & 0 & 2 & 1 & 1 \\
&  &  & 2 & 1448.4 & 579.4 & 413.8 & 60.98 & 14.48 & 46.50 & 0 & 0 & 0 & 0 & 2 & 0 & 1 \\
				
			\hline			
	\end{tabular}}
	\caption{Details of the non-dominated paths found for the instances of 1500 km with w}
	\label{tab:results1500_3}
\end{table}

\endgroup

The results obtained with the algorithm mimicking the real-life behavior of the truck drivers show that the solution found is generally dominated by all the solutions found by the proposed algorithm. A detailed comparison of the solutions found by the two algorithms is presented in Table \ref{tab:comparison}. For each instance, the table reports the maximum savings obtained by the proposed algorithm with respect to the current practice both in terms of absolute value and percentage. Furthermore, the ratio of solutions that dominate the current practice, in both criteria, is reported. For such solutions, the four rightmost columns report the average savings, both in terms of absolute value and percentage over the current practice. It can be observed that the proposed algorithm provides fuel cost savings up to around 115 \euro{} (instance 1500\_b\_10) and a time saving of up to one hour and 30 minutes (instance 1500\_b\_50). In terms of percent deviations, we observe fuel cost savings up to 30\% (e.g., instance 500\_b\_10) and time savings of up to 12\% (e.g., instance 500\_d\_10). Thus, benefits in some instances are indeed significant, especially in terms of fuel cost savings. It is also important to note that another valuable benefit provided by the algorithm is the fact that it generally produces different non-dominated solutions. This is especially valuable when carriers may be interested to evaluate not only the cheapest option (in terms of fuel cost) but other solutions, leading to possibly faster routes.

A note must be made on the results found in the 1000 km instances with $f_o \in \{75, 100\}$. In particular, it can be observed that in instances 1000\_b\_75 and 1000\_b\_100, the solutions found by the two algorithms are practically equivalent, with the difference that, because of its slightly more conservative nature, the proposed algorithm finds different rest locations that are slightly further away from the highway than those found by the CP algorithm. In the remaining instances, i.e., those of the ``d'' and ``w'' scenarios, the same explanation applies for the fuel cost of the solutions, while the substantial time saving of the solutions found by the presented algorithm is due to its ability to smartly combine different kinds of rest stops.

Given the nature of the generated instances, we see fit to report a few statistical measures related to the savings reported in Table \ref{tab:comparison}. Since the initial fastest paths have different lengths, we have scaled the savings reported in Table \ref{tab:comparison} to represent averages over 500 km. Considering these scaled savings over the 122 non-dominated paths (found in the 45 instances), the average fuel savings and time savings are 17.7 \euro{} and 0.08 hours. In this case, the 95\% confidence intervals of the average fuel savings and time savings are $[14.9,20.6]$ and $[-0.06,0.22]$. To further explore these results, we performed a paired t-test (with a 0.05 significance level) to compare the mean difference between the fuel costs of the non-dominated paths and those of the current practice. We conclude that the difference between the two means is statistically significant. We have also performed a similar test between the duration of the non-dominated paths and those of the current practice. We conclude that there is no significant difference in the means of these samples.

\begin{table}[ht]
    \centering
    \resizebox*{!}{\textheight}{%
		\begin{tabular}{rccrrrrrrrrr}
		\multicolumn{3}{c}{Instance} &
		\multicolumn{2}{c}{\begin{tabular}[r]{@{}c@{}}Max absolute \\ saving \end{tabular}} &
		\multicolumn{2}{c}{\begin{tabular}[r]{@{}c@{}}Max percentage \\ saving \end{tabular}} &
		\multicolumn{1}{c}{}
		&
		\multicolumn{2}{c}{\begin{tabular}[r]{@{}c@{}}Avg. absolute \\ saving over CP \end{tabular}} &
		\multicolumn{2}{c}{\begin{tabular}[r]{@{}c@{}}Avg. percentage \\ saving over CP \end{tabular}} 
		\\
		\cmidrule(lr){1-3} \cmidrule(lr){4-5} \cmidrule(lr){6-7} \cmidrule(lr){9-10} \cmidrule(lr){11-12}
		km & HOS & fuel & 
		\begin{tabular}[r]{@{}c@{}}Fuel \\ (\euro) \end{tabular} & 
		\begin{tabular}[r]{@{}c@{}}Time \\ (h) \end{tabular} & 
		\begin{tabular}[r]{@{}c@{}}Fuel \\ (\%) \end{tabular} & 
		\begin{tabular}[r]{@{}c@{}}Time \\ (\%) \end{tabular} & 
		\begin{tabular}[r]{@{}c@{}} Routes \\ Dom. \\ CP  \end{tabular} &
		\begin{tabular}[r]{@{}c@{}}Fuel \\ (\euro) \end{tabular} & 
		\begin{tabular}[r]{@{}c@{}}Time \\ (h) \end{tabular} & 
		\begin{tabular}[r]{@{}c@{}}Fuel \\ (\%) \end{tabular} & 
		\begin{tabular}[r]{@{}c@{}}Time \\ (\%) \end{tabular}
		\\
		
		\hline \hline
		
		500 & b & 10 & $-$62.06 & $-$0.73 & $-$29.98 & $-$4.43 & 2/2 & $-$60.33 & $-$0.59 & $-$29.14 & $-$3.61 \\
500 & b & 25 & $-$29.79 & $-$0.47 & $-$14.83 & $-$2.89 & 5/5 & $-$26.08 & $-$0.39 & $-$12.98 & $-$2.41 \\
500 & b & 50 & $-$17.70 & $-$0.44 & $-$8.57 & $-$2.74 & 1/1 & $-$17.70 & $-$0.44 & $-$8.57 & $-$2.74 \\
500 & b & 75 & $-$17.70 & $-$0.44 & $-$8.57 & $-$2.74 & 1/1 & $-$17.70 & $-$0.44 & $-$8.57 & $-$2.74 \\
500 & b & 100 & $-$17.70 & $-$0.44 & $-$8.57 & $-$2.74 & 1/1 & $-$17.70 & $-$0.44 & $-$8.57 & $-$2.74 \\
\hline
500 & d & 10 & $-$62.06 & $-$0.73 & $-$29.98 & $-$11.80 & 2/2 & $-$60.33 & $-$0.59 & $-$29.14 & $-$9.60 \\
500 & d & 25 & $-$29.79 & $-$0.47 & $-$14.83 & $-$7.89 & 5/5 & $-$26.08 & $-$0.39 & $-$12.98 & $-$6.57 \\
500 & d & 50 & $-$17.70 & $-$0.44 & $-$8.57 & $-$7.48 & 1/1 & $-$17.70 & $-$0.44 & $-$8.57 & $-$7.48 \\
500 & d & 75 & $-$17.70 & $-$0.44 & $-$8.57 & $-$7.48 & 1/1 & $-$17.70 & $-$0.44 & $-$8.57 & $-$7.48 \\
500 & d & 100 & $-$17.70 & $-$0.44 & $-$8.57 & $-$7.48 & 1/1 & $-$17.70 & $-$0.44 & $-$8.57 & $-$7.48 \\
\hline
500 & w & 10 & $-$62.06 & $-$0.73 & $-$29.98 & $-$11.80 & 2/2 & $-$60.33 & $-$0.59 & $-$29.14 & $-$9.60 \\
500 & w & 25 & $-$29.79 & $-$0.47 & $-$14.83 & $-$7.89 & 5/5 & $-$26.08 & $-$0.39 & $-$12.98 & $-$6.57 \\
500 & w & 50 & $-$17.70 & $-$0.44 & $-$8.57 & $-$7.48 & 1/1 & $-$17.70 & $-$0.44 & $-$8.57 & $-$7.48 \\
500 & w & 75 & $-$17.70 & $-$0.44 & $-$8.57 & $-$7.48 & 1/1 & $-$17.70 & $-$0.44 & $-$8.57 & $-$7.48 \\
500 & w & 100 & $-$17.70 & $-$0.44 & $-$8.57 & $-$7.48 & 1/1 & $-$17.70 & $-$0.44 & $-$8.57 & $-$7.48 \\
\hline
\hline
1000 & b & 10 & $-$66.13 & $-$0.18 & $-$17.01 & $-$0.86 & 1/2 & $-$61.81 & $-$0.18 & $-$15.90 & $-$0.86 \\
1000 & b & 25 & $-$65.65 & 0.00 & $-$15.95 & $-$0.01 & 1/5 & $-$18.92 & 0.00 & $-$4.60 & $-$0.01 \\
1000 & b & 50 & $-$10.11 & $-$0.01 & $-$2.71 & $-$0.04 & 1/3 & $-$1.96 & $-$0.01 & $-$0.53 & $-$0.04 \\
1000 & b & 75 & 0.77 & 0.02 & 0.21 & 0.09 & 0/1 &  &  &  &  \\
1000 & b & 100 & 0.77 & 0.02 & 0.21 & 0.09 & 0/1 &  &  &  &  \\
\hline
1000 & d & 10 & $-$66.13 & $-$0.93 & $-$17.01 & $-$4.23 & 2/2 & $-$63.97 & $-$0.74 & $-$16.46 & $-$3.37 \\
1000 & d & 25 & $-$65.65 & $-$0.75 & $-$15.95 & $-$3.45 & 5/5 & $-$49.20 & $-$0.68 & $-$11.95 & $-$3.10 \\
1000 & d & 50 & $-$10.11 & $-$0.01 & $-$2.71 & $-$0.04 & 1/3 & $-$1.96 & $-$0.01 & $-$0.53 & $-$0.04 \\
1000 & d & 75 & 0.77 & $-$0.73 & 0.21 & $-$3.36 & 0/1 &  &  &  &  \\
1000 & d & 100 & 0.77 & $-$0.73 & 0.21 & $-$3.36 & 0/1 &  &  &  &  \\
\hline
1000 & w & 10 & $-$66.13 & $-$0.93 & $-$17.01 & $-$1.66 & 2/2 & $-$63.97 & $-$0.74 & $-$16.46 & $-$1.33 \\
1000 & w & 25 & $-$65.65 & $-$0.75 & $-$15.95 & $-$1.35 & 5/5 & $-$49.20 & $-$0.68 & $-$11.95 & $-$1.21 \\
1000 & w & 50 & $-$10.11 & $-$0.01 & $-$2.71 & $-$0.02 & 1/3 & $-$1.96 & $-$0.01 & $-$0.53 & $-$0.02 \\
1000 & w & 75 & 0.77 & $-$0.73 & 0.21 & $-$1.31 & 0/1 &  &  &  &  \\
1000 & w & 100 & 0.77 & $-$0.73 & 0.21 & $-$1.31 & 0/1 &  &  &  &  \\
\hline
\hline
1500 & b & 10 & $-$114.76 & $-$1.36 & $-$17.93 & $-$3.49 & 2/2 & $-$114.20 & $-$1.12 & $-$17.84 & $-$2.88 \\
1500 & b & 25 & $-$75.10 & $-$1.34 & $-$13.41 & $-$3.47 & 6/8 & $-$37.15 & $-$1.05 & $-$6.63 & $-$2.73 \\
1500 & b & 50 & $-$62.07 & $-$1.46 & $-$10.80 & $-$3.79 & 3/3 & $-$48.53 & $-$1.22 & $-$8.44 & $-$3.16 \\
1500 & b & 75 & $-$33.41 & $-$0.59 & $-$5.63 & $-$1.56 & 2/2 & $-$30.89 & $-$0.45 & $-$5.20 & $-$1.18 \\
1500 & b & 100 & $-$24.20 & $-$1.36 & $-$4.01 & $-$3.51 & 1/1 & $-$24.20 & $-$1.36 & $-$4.01 & $-$3.51 \\
\hline
1500 & d & 10 & $-$109.00 & $-$1.22 & $-$17.18 & $-$4.28 & 2/4 & $-$108.04 & $-$0.98 & $-$17.03 & $-$3.43 \\
1500 & d & 25 & $-$71.34 & $-$1.23 & $-$12.82 & $-$4.37 & 6/8 & $-$33.38 & $-$0.95 & $-$6.00 & $-$3.36 \\
1500 & d & 50 & $-$61.57 & $-$0.70 & $-$10.72 & $-$2.54 & 3/4 & $-$46.40 & $-$0.47 & $-$8.08 & $-$1.70 \\
1500 & d & 75 & $-$38.57 & $-$1.22 & $-$6.44 & $-$4.32 & 2/2 & $-$36.04 & $-$1.07 & $-$6.02 & $-$3.81 \\
1500 & d & 100 & $-$19.15 & $-$1.22 & $-$3.20 & $-$4.32 & 1/2 & $-$18.68 & $-$1.22 & $-$3.12 & $-$4.32 \\
\hline
1500 & w & 10 & $-$109.00 & $-$1.22 & $-$17.18 & $-$1.95 & 2/4 & $-$108.04 & $-$0.98 & $-$17.03 & $-$1.56 \\
1500 & w & 25 & $-$71.34 & $-$1.23 & $-$12.82 & $-$1.98 & 6/8 & $-$33.38 & $-$0.95 & $-$6.00 & $-$1.52 \\
1500 & w & 50 & $-$61.57 & $-$0.70 & $-$10.72 & $-$1.14 & 3/4 & $-$46.40 & $-$0.47 & $-$8.08 & $-$0.76 \\
1500 & w & 75 & $-$38.57 & $-$1.22 & $-$6.44 & $-$1.96 & 2/2 & $-$36.04 & $-$1.07 & $-$6.02 & $-$1.73 \\
1500 & w & 100 & $-$19.15 & $-$1.22 & $-$3.20 & $-$1.96 & 1/2 & $-$18.68 & $-$1.22 & $-$3.12 & $-$1.96 \\
		\hline \\
	\end{tabular}}
    \caption{Summary of the savings as differences with the current practice solution}
    \label{tab:comparison}
\end{table}

\begin{figure}[ht]
	\centering
	\includegraphics[scale=1]{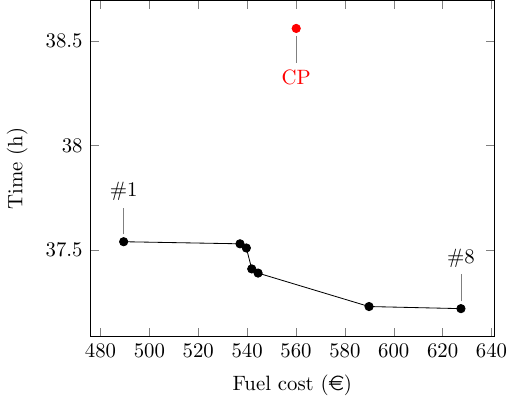}
	\caption{Pareto front and practice solution (red) for instance 1500\_b\_25}
	\label{fig:front}
\end{figure}

\begin{figure}[ht]
	\centering
	\includegraphics[width=\textwidth]{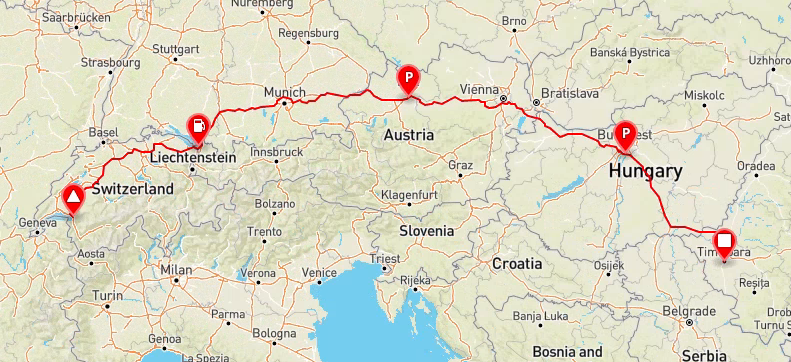}
	\caption{Path \# 1 for the 1500\_b\_25 instance}
	\label{fig:example1}
\end{figure}

\begin{figure}[ht]
	\centering	
	\includegraphics[width=\textwidth]{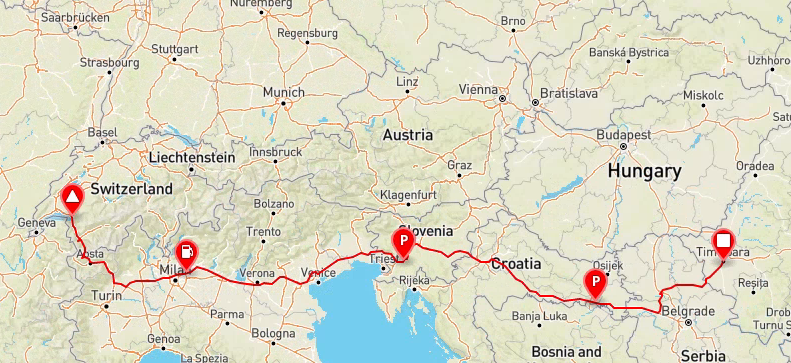}
	\caption{Path \# 8 for the 1500\_b\_25 instance}
	\label{fig:example2}
\end{figure}

\begin{figure}[ht]
	\centering	
	\includegraphics[width=\textwidth]{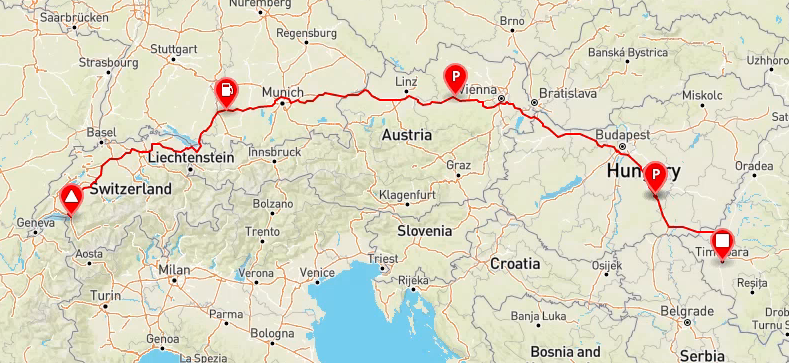}
	\caption{Path CP for the 1500\_b\_25 instance } 
	\label{fig:exampleCP}
\end{figure}

We further illustrate the differences between the obtained paths from the proposed algorithm and the path obtained by current practice through an example. Figure \ref{fig:front} reports the Pareto front of the solutions found by the algorithm for instance 1500\_b\_25 and, in red, the solution found by the current practice. The figure clearly shows the  gap between the solution from practice and all non-dominated solutions found by the algorithm. Two of the obtained paths (\#1 and \#8) are illustrated in Figures \ref{fig:example1} and \ref{fig:example2}. We observe that the paths follow very different trajectories. As observed in Table \ref{tab:results1500_2}, while \#8 is shorter than \#1 by about 20 minutes, because of the higher fuel prices, \#1 is about 139 \euro{} cheaper than \#8. The path obtained by current practice is presented in Figure \ref{fig:exampleCP}. While the CP path follows a road itinerary somewhat similar to that of path \#1, the ability of the proposed algorithm to evaluate stops at different fuel and rest locations leads to cheaper refueling prices and shorter detours from the main highways.


\section{Managerial implications and conclusions }\label{sec:conclusions}
In Section \ref{MI} we draw managerial insights derived from the algorithm, experiments, and discuss the main steps towards a practical implementation of our methodology. In Section \ref{CN} we present the overall conclusions of the paper and discuss future research directions.

\subsection{Managerial implications} \label{MI}

The variability of fuel prices from one fuel station to another entails that gains can be achieved by route deviations. Such gains grow with the increase in the planned route length. Indeed, our results show that the maximum fuel savings achieved with the 1500 km routes is about 1.85 times more than the maximum fuel savings achieved with the 500 km routes. Therefore, accounting for fuel prices is particularly important in the context of long-haul transportation. Fuel savings are mostly observed when deviations from the fastest routes occur, which increase travel times. Therefore, having a bi-objective planning tool is instrumental in accounting for fuel prices. However, we note that the number of non-dominated paths is generally rather limited. In our experiments, the average number of non-dominated paths is 2.7. Such a limited number of paths allows visually inspecting them, and thus further facilitates managerial decisions.

Our algorithm resulted in substantial savings when compared to the current practice of the company involved. Considering all non-dominated paths, the average savings are 35.7~\euro{}, while the average path duration increases by 3.5 minutes. However, as shown in Section \ref{sec:results}, the latter increase is not statistically significant. Moreover, in 91 out of the 122 cases, the non-dominated paths dominate the current practice path on both criteria.

We aimed at developing a prototype tool that is based on free readily available information and interfaces. In particular, our algorithm interacts in a seamless fashion with \gis{} tools, i.e., routing libraries, servers, and web-based data filtering tools. By doing so, we traded simplicity and affordability with speed and location accuracy. Indeed, on average more than 90\% of the runtime is spent retrieving the relevant information from the \gis{} tools. However, we note that in a practical context, if the algorithm is interfaced with commercial and properly configured \gis{} tools, its run times would substantially decrease. For example, simply running our algorithm with an alternative OverPass server (among those provided by the open-source community) almost doubles the runtimes. Another point that would require careful considerations is the search strategy for the stop locations. In particular, the trade-off between a more accurate search strategy yielding fewer stop locations (such as isodistance or isochrone searches), potentially reducing the running time of the algorithm, and a simple and rapid search such as the presented radius search, leaving the assessment of the quality and feasibility of the candidate stop locations to the algorithm, would need to be carefully evaluated. We believe that an isochrone search strategy becomes highly relevant in dynamic settings.

Our prototype tool was designed to work with open-source data. However, acquiring two sources of input data may produce more accurate results, and thus improve practical applications of the prototype tool. First, obtaining rest locations, which are not refueling stations, could enhance the results. Second, procuring daily fuel prices at stations across a number of countries would certainly render the results more applicable. 

\subsection{Conclusions} \label{CN}
Fuel cost is one of the major cost components in long-haul transportation. However, the overwhelming majority of the scientific literature on optimizing long-haul truck transportation ignores refueling decisions. By doing so, this literature implicitly assumes that fuel prices are equal throughout the road network, and thus vehicles should not deviate from their fastest paths. This is a simplistic assumption as fuel prices may vary substantially from one refueling station to another. Therefore, possible detours for refueling operations may yield substantial cost savings. Moreover, such detours may entail a change in the trajectory of the planned path.

In this paper, we studied the long-haul truck scheduling problem on a road network where fuel is accounted for as a cost component and the driver is subject to \hos{} regulations. An origin-destination path has to be determined. This path includes the locations where the driver should stop either for \hos{} compliance, refueling or a combination of the two. Two objectives are considered: the minimization of fuel cost and of path duration.
Given the size of long-haul routes, the number of possible paths, the number of potential refueling stations and the number of potential rest locations, the algorithm we designed interacts with a \gis{} to identify relevant paths and stop locations.
For a set of real-sized instances with distances ranging from 500 to 1500 km, non-dominated paths are generated and compared with the path generated by the current practice of a logistics company.

The problem and the solution approach proposed in this paper concerns \hos{} regulations coming from the European Union legislation. However, the approach can be modified to take into account different regulations associated with other legislations. In fact, it suffices to modify the parameters associated with maximum driving and duty time. Rest locations will be then identified accordingly.

This paper opens a new direction for research in transportation, oriented to considering the refueling cost, that differs from station to station. In particular, given that in some countries fuel prices may change during the day, an extension of our algorithm to handle dynamic prices is a highly relevant research avenue. Considering a dynamic setting also applies to traveling times on the road network. Planning routing and stopping according to traffic conditions could provide substantial benefits.
Another relevant research direction to which this paper contributes is related to the design of optimization algorithms on road networks that interact with a \gis{}.

\section*{Acknowledgements}

The authors would like to thank the anonymous referees who provided useful and detailed comments on a previous version of the manuscript.

\bibliographystyle{abbrvnat}
\bibliography{biblio}

\end{document}